\documentclass[preprint,12pt,3p,fleqn]{elsarticle}
\usepackage{amssymb}
\usepackage{amsthm}
\usepackage{amsmath}
\usepackage{multirow}
\usepackage{booktabs}
\usepackage{graphicx}
\usepackage{epsfig}
\usepackage{epstopdf}
\usepackage{pdflscape}
\usepackage{color}
\usepackage{amsthm}
\usepackage{amsmath}
\usepackage[english]{babel}
\usepackage{subfigure}
\usepackage{threeparttable}
\usepackage[normalem]{ulem}
\usepackage{soul}
\usepackage{cancel}

\usepackage{graphicx}
\usepackage{subfigure}
\usepackage{float}
\usepackage{tabularx}
\usepackage{graphicx}
\usepackage{threeparttable}
\usepackage[figuresright]{rotating}
\allowdisplaybreaks

\theoremstyle{definition}

\theoremstyle{remark}

\journal{Journal of Computational Physics}

\begin{document}
\begin{frontmatter}

\title{A low-dissipation shock-capturing framework with flexible nonlinear dissipation control}
\author[a]{Yue Li}
\ead{yue06.li@tum.de}

\author[b]{Lin Fu\corref{cor1}}
\ead{linfu@stanford.edu}
\cortext[cor1]{Corresponding author.}

\author[a]{Nikolaus A. Adams}
\ead{nikolaus.adams@tum.de}

\address[a]{Chair of Aerodynamics and Fluid Mechanics, Department of Mechanical Engineering, Technical University of Munich, 85748 Garching, Germany}

\address[b]{Center for Turbulence Research, Stanford University, Stanford, CA 94305, USA}

\begin{abstract}
In this work, a framework to construct arbitrarily high-order low-dissipation shock-capturing schemes with flexible and controllable nonlinear dissipation for convection-dominated problems is proposed. While a set of candidate stencils of incremental width is constructed, each one is indicated as smooth or nonsmooth by the ENO-like stencil selection procedure proposed in the targeted essentially non-oscillatory (TENO) scheme [Fu et al., Journal of Computational Physics 305 (2016): 333-359]. Rather than being discarded directly as with TENO schemes, the nonsmooth candidates are filtered by an extra nonlinear limiter, such as a monotonicity-preserving (MP) limiter or a total variation diminishing (TVD) limiter. Consequently,  high-order reconstruction is achieved by assembling candidate fluxes with optimal linear weights since they are either smooth reconstructions or filtered ones which feature good non-oscillation property. A weight renormalization procedure as with the standard TENO paradigm is not necessary. This new framework concatenates the concepts of TENO, WENO and other nonlinear limiters for shock-capturing, and provides a new insight to designing low-dissipation nonlinear schemes. Through the adaptation of nonlinear limiters, nonlinear dissipation in the newly proposed framework can be controlled separately without affecting the performance in smooth regions. Based on the proposed framework, a family of new six- and eight-point nonlinear schemes with controllable dissipation is proposed. A set of critical benchmark cases involving strong discontinuities and broadband fluctuations is simulated. Numerical results reveal that the proposed new schemes capture discontinuities sharply and resolve the high-wavenumber fluctuations with low dissipation, while maintaining the desired accuracy order in smooth regions.
\end{abstract}


\end{frontmatter}

% \linenumbers

\section{Introduction}
High-order and high-resolution shock-capturing schemes are popular numerical methods to solve compressible fluid problems, which may involve discontinuities. However, there still is a need for flexible environments to design such schemes that are highly accurate in smooth regions with low numerical dissipation and can capture discontinuities with low parasitic noise in smooth regions near discontinuous.

Among all the schemes, weighted essentially non-oscillatory (WENO) schemes, first proposed by Liu et al. \cite{Liu1994} and further improved by Jiang and Shu \cite{Jiang1996}, probably are the most  widely accepted discretization schemes. These schemes are developed from the essentially non-oscillatory (ENO) concept \cite{harten1987uniformly}, which selects the smoothest stencil from a set of candidates to prevent the reconstruction from crossing discontinuities. Instead of selecting the smoothest candidate, WENO deploys a convex combination of all candidate stencils to achieve high-order accuracy in smooth regions. In nonsmooth regions, reduced weights are assigned to the stencils through nonlinear smoothness indicators. Thus, spurious numerical oscillations can be largely avoided. However, practical implementations reveal that the WENO-JS scheme \cite{Jiang1996} is rather dissipative and may lose the accuracy order near critical points due to the strong nonlinear adaptation within the weighting strategy.

In order to reduce dissipation while maintaining high-order accuracy in smooth regions and shock-capturing capability, two main methods have been introduced: (i) decreasing the nonlinear adaptive dissipation \cite{henrick2005mapped}\cite{borges2008improved}\cite{hill2004hybrid}, and (ii) optimizing the background counterpart linear scheme \cite{martin2006bandwidth}\cite{Hu2010}. Henrick et al. \cite{henrick2005mapped} suggest that restoring the optimal accuracy order near critical points can resolve the over-dissipation issue. Necessary and sufficient conditions have been derived and a mapping strategy based on classical WENO, referred to as WENO-M, has been developed accordingly. However, as investigated by Borges et al. \cite{borges2008improved}, the improvements demonstrated by WENO-M \cite{henrick2005mapped} over WENO-JS result from larger weights of nonsmooth candidate stencils and less nonlinear adaptive dissipation, rather than from higher accuracy near critical points. To further decrease the nonlinear adaptive dissipation, a new smoothness indicator has been proposed by Borges et al. \cite{borges2008improved}, which introduces a global high-order undivided difference into the weighting strategy of the WENO-JS scheme. As an alternative approach \cite{hill2004hybrid}, the nonlinear dissipation adaptation can be controlled by comparing the ratio between the largest and the smallest calculated smoothness indicator. By introducing a problem-dependent threshold, adaptation is eliminated when the ratio is below the threshold. The aforementioned WENO-like schemes reduce the dissipation of the classical WENO significantly. However, they are still too dissipative for resolving turbulence-like high-wavenumber fluctuations. Besides decreasing the nonlinear adaptive dissipation, optimizing the background counterpart linear scheme of WENO also can improve the overall dissipation property. The WENO-SYMOO \cite{martin2006bandwidth} and WENO-CU6 \cite{Hu2010} schemes optimize the background scheme as a central scheme by introducing the contribution of the downwind stencil. However, anti-dissipation is inherently built in for a certain wavenumber range, which may be problematic for critical applications. Therefore, in terms of designing low-dissipation shock capturing schemes, the strategy of choosing a stable linear scheme with upwind-biased candidate stencils as a background scheme and then optimizing the nonlinear dissipation adaptation to its limit is preferable.

To further reduce the dissipation by suppressing nonlinear adaptation and achieve stable shock-capturing capability, a family of high-order targeted ENO (TENO) schemes has been proposed by Fu et al. \cite{fu2016family}\cite{fu2017targeted}\cite{fu2018new}\cite{fu2019low}. Unlike the WENO-like smooth convex combination, TENO either adopts a candidate stencil with its optimal weight for the final reconstruction or discards it completely when crossed by a genuine discontinuity. This procedure is ensured by a cut-off parameter $C_T$, which identifies the smooth and nonsmooth regions. In smooth regions, nonlinear adaptation is unnecessary and the numerical dissipation is reduced to that of the linear background scheme. In nonsmooth regions, nonlinear dissipation is introduced to suppress oscillations and to achieve the robust shock-capturing capability. This strategy successfully reduces the numerical dissipation so that turbulence-like high-wavenumber fluctuations can be sustained \cite{fu2017targeted}. However, based on many observations \cite{fu2017targeted}\cite{fu2019low}, tailored nonlinear dissipation instead of constant low dissipation is needed in nonsmooth regions to mitigate numerical artifacts caused by the low-dissipation linear background schemes.

{
An alternative solution is to introduce the hybrid concept, for which the efficient linear scheme is applied in smooth regions while the expensive nonlinear scheme is adopted in the vicinity of discontinuities. Adams and Shariff \cite{adams1996high} have introduced the concept of hybridization of a high-resolution compact scheme with the nonlinear high-order ENO scheme. Pirozzoli \cite{pirozzoli2002conservative} proposes a conservative hybrid scheme for resolving the compressible turbulent flows by combining a compact scheme with WENO. However, the performance of these hybrid methods strongly depends on an effective discontinuity indicator \cite{hu2015efficient}, which is mostly case-sensitive.}

In this paper, we propose an extended framework based on the TENO concept. The main objective is to achieve adaptive nonlinear dissipation property without deteriorating the spectral resolution properties of TENO in low-wavenumber regions. The key idea is that rather than discarding completely the nonsmooth candidate stencils, their contributions are filtered by a nonlinear limiter. The filtered candidate stencil contributions are assembled to form the final reconstruction with their respective optimal linear weights. On smooth candidate stencils, the nonlinear limiter will not be activated and consequently the properties of TENO schemes in smooth regions are maintained. Different choices of the nonlinear limiter are possible, such that a new family of schemes with adjustable nonlinear dissipation in nonsmooth regions is obtained.

The rest of this paper is organized as follows. In section 2, the concept of the finite-difference high-order reconstruction schemes and the construction of standard TENO schemes are briefly reviewed. In section 3, a general framework to construct nonlinear shock-capturing schemes is proposed and the filtering strategy for nonsmooth candidate stencils are elaborated. In section 4, a set of benchmark cases is considered for assessment of the developed schemes. The concluding remarks are given in the last section.

\section{Concepts of standard TENO schemes for hyperbolic conservation laws}

To facilitate the presentation, we consider a one-dimensional scalar hyperbolic conservation law
\begin{equation}
\label{eq:conservationlaw}
\frac{\partial u }{\partial t} + \frac{\partial}{\partial x}f(u) = 0 ,
\end{equation}
where $u$ and $f$ denote the conservative variable and the flux function respectively. The characteristic signal speed is assumed to be positive $\frac{\partial f(u)}{\partial u} > 0$. The developed schemes can be extended to systems of conservation laws and multi-dimensional problems in a straight-forward manner.

For a uniform Cartesian mesh with cell centers ${x_i} = i\Delta x$ and cell interfaces ${x_{i+1/2}} = {x_i} + \frac{\Delta x}{2}$, the spatial discretization results in a set of ordinary differential equations
\begin{equation}
\label{eq:ordinaryEq}
\frac{d{u_i(t)}}{dt} =  - \frac{\partial f}{\partial x}\left| {_{x = {x_i}}} \right., \text{ } i = 0,\cdots,n,
\end{equation}
where ${u_i}$ is a numerical approximation to the point value $u(x_i,t)$. The semi-discretization with  Eq.~(\ref{eq:ordinaryEq}) can be further discretized by a conservative finite-difference scheme as
\begin{equation}
\label{eq:implicitEq}
\frac{d{u_i}}{dt} =  - \frac{1}{\Delta x}({h_{i + 1/2}} - {h_{i - 1/2}}) ,
\end{equation}
where the primitive function $h(x)$ is implicitly defined by
\begin{equation}
\label{eq:definitionEq}
f(x) = \frac{1}{\Delta x}\int_{x - \Delta x/2}^{x + \Delta x/2} {h(\xi )d\xi },
\end{equation}
and ${h_{i \pm 1/2}} = h({x_i} \pm \frac{{\Delta x}}{2})$. A high-order approximation of $h(x)$ at the cell interface has to be reconstructed from the cell-averaged values of $f(x)$ at the cell centers. Eq.~(\ref{eq:implicitEq}) can be written as
\begin{equation}
\label{eq:approximateEq}
\frac{d{u_i}}{dt} \approx  - \frac{1}{\Delta x}({\widehat f_{i + 1/2}} - {\widehat f_{i - 1/2}}) ,
\end{equation}
where ${\hat f_{i \pm 1/2}}$ denotes the high-order approximate numerical fluxes and can be computed from a convex combination of $K - 2$ candidate-stencil fluxes
\begin{equation}
\label{eq:convex}
\widehat f_{i + 1/2} = \sum\limits_{k = 0}^{K - 3} {w_k} \widehat f_{k,i + 1/2} .
\end{equation}

A $({r_k} - 1)$-degree polynomial is assumed for each candidate stencil as
\begin{equation}
\label{eq:approximatepolynomial}
h(x) \approx {\hat f_k}(x) = \sum\limits_{l = 0}^{{r_k} - 1} {{a_{l,k}}} {x^l},
\end{equation}
where ${r_k}$ denotes the point number of candidate stencil $k$. After substituting Eq.~(\ref{eq:approximatepolynomial}) into Eq.~(\ref{eq:definitionEq}) and evaluating the integral functions at the stencil nodes, the coefficients ${a_{l,k}}$ are determined by solving the resulting system of linear algebraic equations.

For hyperbolic conservation laws, discontinuities may occur even when the initial condition is smooth enough.
The challenge is to develop a reconstruction scheme which is high-order accurate in smooth regions and captures discontinuities sharply and stably in nonsmooth regions. In the following, we recall essential elements of the recently proposed TENO \cite{fu2016family}\cite{fu2017targeted}\cite{fu2018new}\cite{fu2019low}\cite{Fu2019hybrid} concept.
\subsection{Candidate stencil arrangement}

\begin{figure}[htbp]
\centering
\subfigure{\includegraphics[width=0.8\textwidth]{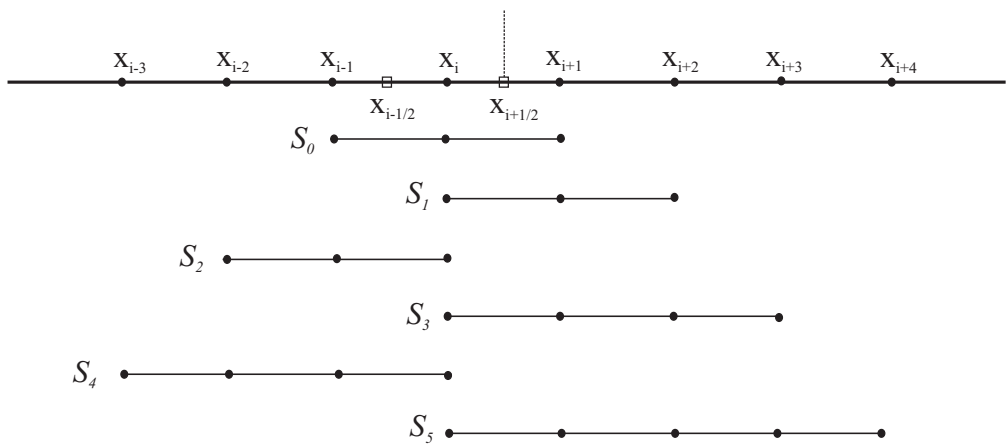}}
\caption{Candidate stencils with incremental width towards high-order reconstruction. While the first three stencils are identical to that of the classical fifth-order WENO-JS scheme, all candidate stencils possess at least one upwind point. The maximum $K$th-order scheme can be constructed by combining candidates ranging from $S_0$ to $S_{K-3}$. The advantage of such candidate stencil arrangement over that of classical WENO schemes is the robustness improvement in very-high-order versions.}
\label{Fig:incremental_stencil}
\end{figure}

As shown in Fig.~\ref{Fig:incremental_stencil}, arbitrarily high-order TENO schemes are constructed from a set of candidate stencils with incremental width \cite{fu2016family}.
The sequence of stencil width $r$ varying versus the accuracy order $K$ is as
\begin{equation}
\label{eq:rdistribution}
\{ {r_k}\}  = \left\{ {\begin{array}{*{20}{c}}
{\{ \underbrace {3,3,3,4,\cdots,\frac{{K + 2}}{2}}_{0,\cdots,K - 3}\} ,}&{\text{if } \bmod (K,2) = 0,}\\
{\{ \underbrace {3,3,3,4,\cdots,\frac{{K + 1}}{2}}_{0,\cdots,K - 3}\} ,}&{\text{if } \bmod (K,2) = 1.}
\end{array}} \right.
\end{equation}

Note that, although the general candidate stencil arrangement of classical WENO schemes has some limitations \cite{fu2016family}, the following ideas discussed will still be applicable even if it is employed.

\subsection{Scale-separation procedure}

To isolate effectively discontinuities from smooth regions, smoothness indicators with strong scale-separation capability are given as \cite{fu2016family}
\begin{equation}
\label{eq:measureTENO}
{\gamma _k} = {\left( C + \frac{{{\tau _K}}}{{{\beta _{k,r}} + \varepsilon }}\right)^q}{\text{ , }}k = 0,\cdots,K - 3 ,
\end{equation}
where $\varepsilon  = {10^{ - 40}}$ is introduced to prevent the zero denominator. The parameters $C = 1$ and $q = 6$ are chosen for strong scale separation. Following Jiang and Shu \cite{Jiang1996}, by counting all the possible high-order undivided differences, $\beta _{k,r}$ can be given as
\begin{equation}
\label{eq:measureTENO}
{\beta _{k,r}} = \sum\limits_{j = 1}^{r - 1} {\Delta {x^{2j - 1}}} \int_{{x_{i - 1/2}}}^{{x_{i + 1/2}}} {{{\left(\frac{{{d^j}}}{{d{x^j}}}{{\hat f}_k}(x)\right)}^2}dx} .
\end{equation}
A sixth-order $\tau _K$, which allows for good stability with a reasonable large CFL number, can be constructed as \cite{fu2017targeted}
\begin{equation}
\label{eq:taukTENO}
\tau _K = \left | {\beta _K} - \frac{1}{6}({\beta _{1,3}} + {\beta _{2,3}} + 4{\beta _{0,3}}) \right | = O(\Delta {x^6}),
\end{equation}
where the ${\beta _K}$ measures the global smoothness on the $K$-point full stencil.
\subsection{ENO-like stencil selection}

For TENO schemes \cite{fu2016family}, the measured smoothness indicators are first normalized as
\begin{equation}
\label{eq:normalize}
{{\chi}_k } = \frac{{{\gamma _k}}}{{\sum\nolimits_{k = 0}^{K - 3} {{\gamma _k}} }} ,
\end{equation}
and subsequently filtered by a sharp cut-off function
\begin{equation}
\label{eq:cut-off}
{\delta _k} = \left\{ {\begin{array}{*{20}{c}}
0, &{\text{if }{\chi}_k  < {C_T},}\\
1, &{\text{otherwise}.}
\end{array}} \right.
\end{equation}
In such a way, all candidate stencils are then identified to be either sufficiently smooth or nonsmooth with a discontinuity crossing the stencil.

\subsection{Nonlinear adaptation strategy for $C_T$}
Although the above weighting strategy is sufficient to separate smooth regions from discontinuities, undesirable numerical dissipation is generated for turbulence-like high-wavenumber fluctuations since they are treated in a similar manner as discontinuities. Previous research reveals that the low-order undivided difference is more sensitive to distinguish high-wavenumber fluctuations from genuine discontinuities than high-order undivided differences \cite{fu2018new}\cite{fu2018targeted}\cite{fu2018improved}.

Motivated by Ren et al.\cite{ren2003characteristic}, the local smoothness of flow field can be indicated by
\begin{equation}
\label{eq:shock_indicator}
\left\{ {\begin{array}{*{20}{c}}
m = 1-{\text{min}(1, \frac{{\eta}_{i+1/2}}{C_r}),}\\
{\eta}_{i+1/2}={\text{min}({\eta}_{i-1},{\eta}_{i},{\eta}_{i+1},{\eta}_{i+2}).}
\end{array}} \right.
\end{equation}
where

\begin{equation}
{{\eta}_i } = \frac{|2\Delta f_{i+1/2}\Delta f_{i-1/2}|+\epsilon}{(\Delta f_{i+1/2})^2+(\Delta f_{i-1/2})^2+ \epsilon} ,
\end{equation}
\begin{equation}
{\Delta f_{i+1/2} } = f_{i+1}-f_{i}  , \\
{\epsilon} = \frac{0.9{C_r}}{1-0.9{C_r}}{\xi}^2 ,
\end{equation}
and the parameters $\xi=10^{-3}$, $C_r = 0.23$ are used for the eight-point TENO schemes \cite{fu2018targeted}.
Low numerical dissipation for turbulence-like high-wavenumber fluctuations is achieved by adjusting the cut-off parameter $C_T$ as
\begin{equation}
\label{eq:cut-off_adjust}
\left\{ {\begin{array}{*{20}{c}}
g(m) = (1-m)^4(1+4m),\\
{\beta}={{\alpha}_1-{\alpha}_2(1-g(m)),}\\
C_T=10^{-[\beta]},
\end{array}} \right.
\end{equation}
where $[\beta]$ denotes the maximum integer which is not larger than $\beta$. $g(m)$ is a smoothing kernel based mapping function, and ${\alpha}_1 = 10.5$, ${\alpha}_2 = 3.5$ are set for the eight-point TENO schemes \cite{fu2018targeted}. The eight-point TENO scheme with $C_T$ adaptation is referred to as TENO8A in this paper.

\subsection{The final high-order reconstruction}
In order to remove contributions from candidate stencils containing discontinuities, optimal weights $d_k$ subjected to the cut-off $\delta _k$ are renormalized as
\begin{equation}
\label{eq:weight_renormalize}
{w_k} = \frac{{d_k}{\delta _k}}{{\sum\nolimits_{k = 0}^{K - 3} {d_k}{\delta _k} }}.
\end{equation}
The $K$th-order reconstructed numerical flux evaluated at cell face $i + \frac{1}{2}$ is assembled as
\begin{equation}
\label{eq:reconstuction}
\hat f_{i + 1/2}^K = \sum\limits_{k = 0}^{K - 3} {w_k} {\hat f_{k,i + 1/2}} .
\end{equation}
\section{New framework with flexible nonlinear dissipation control}
In this section, we describe the details of the new reconstruction framework for TENO with flexible nonlinear dissipation control. The new schemes are referred to as TENO-M schemes hereafter. The advantages of TENO-M schemes are not only in keeping the contribution of nonsmooth stencils but also maintaining the performance improvement of TENO. The possibility of keeping nonsmooth stencils rather than abandoning them entirely is achieved by a flexible filtering procedure. As the nonsmooth stencils are filtered to be oscillation free by an extra limiter, they are able to contribute to the final reconstruction without detriment to computational stability. In regions of smooth stencils, nonlinear limiters will not be activated and consequently the performance of TENO in terms of restoring the background linear scheme in low-wavenumber regions is maintained.

\subsection{Evaluation of the candidate numerical fluxes and label stencils based on smoothness indicator}
As shown in Fig.~\ref{Fig:incremental_stencil}, arbitrarily high-order TENO schemes, i.e. both odd and even order, are constructed from a set of candidate stencils with incremental width \cite{fu2016family}. The numerical flux $\widehat f_{k,i + 1/2} $ of the candidate stencil $S_k$ is evaluated by solving the linear system results from Eq.~(\ref{eq:approximatepolynomial}) and Eq.~(\ref{eq:definitionEq}). The corresponding coefficients for TENO schemes up to eighth-order accuracy are presented in Table.~\ref{tab:stencil_coefficient}.

\begin{table}[htbp]\footnotesize
 \caption{Numerical flux ${\widehat f_{k,i + 1/2}} = \sum\nolimits_m {{c_m}} {f_m}$ evaluated from the candidate stencil $S_k$}
 \centering
\label{tab:stencil_coefficient}
 \begin{tabular*}{0.93\textwidth}{@{\extracolsep{\fill}}ccccccccc}
  \toprule
   \centering
    ${\widehat f_{k,i + 1/2}}$&  ${c_{i-3}}$& ${c_{i-2}}$& ${c_{i-1}}$& ${c_{i}}$& ${c_{i+1}}$& ${c_{i+2}}$& ${c_{i+3}}$& ${c_{i+4}}$\\
  \midrule
   \centering
    $k=0$ & & & ${-}\frac{1}{{6}}$& $\frac{5}{{6}}$& $\frac{2}{{6}}$& & &\\
    \\
    $k=1$ & & & & $\frac{2}{{6}}$& $\frac{5}{{6}}$& ${-}\frac{1}{{6}}$& &\\
    \\
    $k=2$ & & $\frac{2}{{6}}$& ${-}\frac{7}{{6}}$& $\frac{11}{{6}}$& & & &\\
    \\
    $k=3$ & & & & $\frac{3}{{12}}$& $\frac{13}{{12}}$& ${-}\frac{5}{{12}}$& $\frac{1}{{12}}$&\\
    \\
    $k=4$ & ${-}\frac{3}{{12}}$& $\frac{13}{{12}}$& ${-}\frac{23}{{12}}$& $\frac{25}{{12}}$& & & &\\
    \\
    $k=5$ & & & & $\frac{12}{{60}}$& $\frac{77}{{60}}$& ${-}\frac{43}{{60}}$& $\frac{17}{{60}}$& ${-}\frac{3}{{60}}$\\
  \bottomrule
 \end{tabular*}
\end{table}
Based on the smoothness indicator, each candidate stencil is identified as smooth or nonsmooth by the ENO-like stencil selection procedure, i.e. Eq.~(\ref{eq:normalize}) and Eq.~(\ref{eq:cut-off}).

\subsection{Filtering of the nonsmooth stencils with an extra limiter}

{Different from the standard TENO schemes, where the detected nonsmooth stencils will be discarded completely, in the present framework, they are filtered by a specific nonlinear limiter and then adopted for the final reconstruction by:}
\begin{equation}
\label{eq:mp_flux}
{\hat f^{\rm{M}}}_{k,i + 1/2}  =  \left\{ {\begin{array}{*{20}{c}}
{f_{k,i + 1/2}^{\rm{lim}} }&{\text{if $\delta _k=0$, } }\\
{{{\hat f}_{k,i + 1/2}},}&{\rm{otherwise,}}
\end{array}} \right.
\end{equation}
where ${f_{k,i + 1/2}^{\rm{lim}}}$ denotes the numerical flux filtered by a nonlinear limiter. A monotonicity-preserving limiter \cite{suresh1997accurate},  a fifth-order TVD  limiter \cite{kim2005accurate} and a second-order Van Albada limiter \cite{harten1983high}\cite{sweby1984high} are deployed to derive the ${f_{k,i + 1/2}^{\rm{lim}}}$ and are referred to as $f_{k,i + 1/2}^{\rm{MP}}$, $f_{k,i + 1/2}^{\rm{TVD5}}$, and $f_{k,i + 1/2}^{\rm{VA}}$ respectively. Consequently, the nonlinear numerical dissipation property of TENO schemes can be tailored by deploying different nonlinear limiters to the nonsmooth stencils. The accuracy order and low-dissipation property of the resulting TENO schemes are maintained in smooth and low-wavenumber regions.

This new framework connects the concepts of TENO and other nonlinear limiters for shock-capturing. Although only three nonlinear limiters are considered in this paper, any other established limiters can be introduced into Eq.~(\ref{eq:mp_flux}).
\subsubsection{TVD nonlinear limiter}
{Following \cite{van1979towards}, TVD nonlinear limiters enforce the TVD property by introducing the slope function $\phi(r)$ to limit the gradient variation. Several second- and third-order TVD schemes with different slope functions have been developed, e.g. with the Minmod limiter, the Superbee limiter and the Van Leer limiter \cite{zhang2016refined}\cite{kemm2011comparative}\cite{arora1997well}. Here, we take a Van Albada limiter as an example to exemplify its implementation.
Following \cite{van1997comparative}, the slope function of a  Van Albada limiter $\phi(r)_{\rm{VA}}$ is defined as
\begin{equation}
\label{eq:van_Albada_limiter}
\phi_{\rm{VA}}(r) = \frac{2r}{r^2 + 1},
\end{equation}
and the slope ratio $r$ based on a three-cell stencil centered at $i$ is
\begin{equation}
\label{eq:slop_ratio}
% \delta^- = f_i - f_{i-1} , \\ \delta^+ = f_{i+1} - f_{i} ,
r=\delta^+/\delta^-,
\end{equation}
where $\delta^- = f_i - f_{i-1} , \delta^+ = f_{i+1} - f_{i}$.
Therefore, the reconstruction scheme filtered by a Van Albada limiter can be written as
\begin{equation}
\label{eq:vAflux}
f_{i+\frac{1}{2}}^{\rm{VA}} = f_i + \frac{\phi_{\rm{VA}}}{4}[(1-\kappa\phi_{\rm{VA}})\delta^- + (1+\kappa\phi_{\rm{VA}})\delta^+].
\end{equation}

}

% Following \cite{van1979towards}, slope limiters are constructed as functions of slope ratio and the slope ratio based on a three-cells stencil centered at $i$ can be described as follows:
%
% \begin{equation}
% \label{eq:muscl}
% f_{i+\frac{1}{2}} ^{\rm{MUSCL}} = f_i + \frac{1}{4}[(1-\kappa)\phi(r) + (1 + \kappa) r \phi(\frac{1}{r})] (f_{i}-f_{i-1}) ,
% \end{equation}
% \begin{equation}
% \label{eq:slop_ratio}
% % \delta^- = f_i - f_{i-1} , \\ \delta^+ = f_{i+1} - f_{i} ,
% r=\delta^+/\delta^-,
% \end{equation}
%%
% The TVD property is enforced by introducing the slope function $\phi(r)$ to limit the gradient variation. Different $\phi(r)$ functions are defined to provide different numerical dissipation,
% With the Van Albada limiter \cite{van1997comparative},  $\phi(r)$ is defined as

% \begin{equation}
% \label{eq:van_Albada_limiter}
% \phi_{\rm{VA}}(r) = \frac{2r}{r^2 + 1},
% \end{equation}

In smooth regions, a third-order reconstruction is restored by setting $\kappa = 1/3$ and the resulting TENO schemes are referred to as TENO-M-VA.

{However, these reconstructions fail to incorporate concrete information of the variation of $\delta^+$ except for the monotonicity constraint. Kim et al. \cite{kim2005accurate} propose to incorporate the curvature information of $\delta^+$ into TVD limiters by the proper choice of an optimal variation function $\beta$ and develop a class of high-order TVD limiters.  As suggested in \cite{kim2005accurate}, the slope function of a high-order TVD limiter $\phi(r)_{\rm{TVD5}}$ based on a five-cell stencil centered at $i$ is defined as }
\begin{equation}
\label{eq:high order TVD}
\phi_{\rm{TVD5}}(r) = \text{max}(0, \text{min}(\alpha, \alpha r, \beta)) ,
\end{equation}
where $\alpha = 2$ and $\beta$ is given as
\begin{equation}
\label{eq:TVDhighorder2}
\beta = \frac{-2/r_{j-1} + 11 +24 r_{j}- 3 r_j r_{j+1}}{30}.
\end{equation}
A reconstruction scheme filtered by a high-order TVD limiter can be written as
\begin{equation}
\label{eq:TVDbeta}
f_{i+\frac{1}{2}}^{\rm{TVD5}} = f_i + \frac{1}{2} \phi_{\rm{TVD5}}\delta^- .
\end{equation}

The resulting TENO schemes with this fifth-order limiter are referred to as TENO-M-TVD5.

\subsubsection{Monotonicity-preserving limiter}
Suresh and Huynh \cite{suresh1997accurate} propose a monotonicity-preserving method to bound the high-order reconstructed data at cell interface by distinguishing smooth local extrema from genuine discontinuity. The resulting monotonicity-preserving schemes allow for the local extremum to develop in the evaluation of cell interface data and is robust for shock-dominated flows.
{
The minmod function with two arguments is
\begin{equation}
\label{eq:minmod222}
{\rm{minmod}}(x,y) = \frac{1}{2}[{\rm{sgn}}(x) + {\rm{sgn}}(y)] \min(\left| {x} \right|, \left| {y} \right|),
\end{equation}
and the minmod function with four arguments is
\begin{equation}
\label{eq:minmod444}
{\rm{minmod}}(a,b,c,d) = \frac{1}{8}[{\rm{sgn}}(a) + {\rm{sgn}}(b)]\left| [{\rm{sgn}}(a) + {\rm{sgn}}(c)][{\rm{sgn}}(a) + {\rm{sgn}}(d)]\right| \min(\left| {a} \right|, \left| {b} \right|,\left| {c} \right|, \left| {d} \right|).
\end{equation}
}
The median function is
\begin{equation}
\label{eq:median}
{\rm{median}}(x,y,z) = x + {\rm{minmod}}(y - x,z - x) .
\end{equation}
While the curvature at the cell center $i$ can be approximated by
\begin{equation}
\label{eq:curvature}
{d_i} = {f_{i + 1}} - 2{f_i} + {f_{i - 1}} ,
\end{equation}
the curvature measurement at the cell interface $i+ 1/2$ can be defined as
\begin{equation}
\label{eq:maximal_curvature}
{d_{i + 1/2}^{\text{M4}}} = {\rm{minmod}} (4{d_i} - {d_{i + 1}},4{d_{i + 1}} - {d_i},{d_i},{d_{i + 1}}).
\end{equation}
This definition is more restrictive since the room for local extrema to develop is reduced when the ratio $d_{i+1}/{d_i}$ is smaller than $1/4$ or larger than $4$.

In order to define the minimum and maximum bounds of data at the interface $x_{i+1/2}$, the left-side upper limiter is given as
\begin{equation}
\label{eq:upper_limiter}
f_{i + 1/2}^{\text{UL}} = {f_i} + \alpha ({f_i} - {f_{i - 1}}),
\end{equation}
where the choice of $\alpha$ in principle should satisfy the condition that $\rm{CFL} \leq 1/(1+\alpha)$ for stability. In our simulations, $\alpha = 1.25$ is employed to enable a choice of $\rm{CFL} = 0.4$.

The median value of the solution at $x_{i+1/2}$ is given by
\begin{equation}
\label{eq:median_value}
f_{i + 1/2}^{\text{MD}} = \frac{1}{2}({f_i} + {f_{i + 1}}) - \frac{1}{2}d_{i + 1/2}^{\text{MD}} .
\end{equation}

The left-side value allowing for a large curvature in the solution at $x_{i+1/2}$ can be given by
\begin{equation}
\label{eq:large_curvature}
f_{i + 1/2}^{\text{LC}} = {f_i} + \frac{1}{2}({f_i} - {f_{i - 1}}) + \frac{\beta }{3}d_{i - 1/2}^{\text{LC}} ,
\end{equation}
where it is recommended to set $\beta = 4$. Following \cite{suresh1997accurate} and \cite{balsara2000monotonicity}, $d_{i + 1/2}^{\text{MD}} = d_{i + 1/2}^{\text{LC}} = {d_{i + 1/2}^{\text{M4}}}$ is adopted. In numerical validations, the parameters in MP limiters are set as $d_{i + 1/2} = d_{i + 1/2}^{\text{M4}}$ and $\beta = 4$.
%%Note that the parameter $\beta$ and curvature measurement method $d_{i + 1/2}$ will influence the inherent dissipation of the MP limiter.
The bounds are given by
\begin{equation}
\label{eq:value_min_max}
\begin{array}{l}
f_{i + 1/2}^{{\rm{L}}{\rm{,min}}} = \max [\min ({f_i},{f_{i + 1}},f_{i + 1/2}^{{\rm{MD}}}),\min ({f_i},f_{i + 1/2}^{{\rm{UL}}},f_{i + 1/2}^{{\rm{LC}}})],\\
f_{i + 1/2}^{{\rm{L}}{\rm{,max}}} = \min [\max ({f_i},{f_{i + 1}},f_{i + 1/2}^{{\rm{MD}}}),\max ({f_i},f_{i + 1/2}^{{\rm{UL}}},f_{i + 1/2}^{{\rm{LC}}})],
\end{array}
\end{equation}
and the monotonicity preserving value for data at the interface $x_{i+1/2}$ is obtained by limiting the predicted value from other reconstructions as
\begin{equation}
\label{eq:value_monotonicity_preserving}
f_{i + 1/2}^{\rm{MP}} = {\rm{median}}(\widehat{f}_{i + 1/2}, f_{i + 1/2}^{{\rm{L}}{\rm{,min}}}, f_{i + 1/2}^{{\rm{L}}{\rm{,max}}}).
\end{equation}

This monotonicity-preserving limiter is implemented in the new framework through Eq.~(\ref{eq:mp_flux}), and the resulting scheme is referred to as TENO-M-MP.
\subsubsection{Implementation discussions}
{Within the TENO-M framework, the implementation of the MP limiter is slightly different from the TVD-type limiter. With a MP limiter, the minimum and maximum bounds of flux at the interface $x_{i+1/2}$ are first computed. Afterwards, these bounds are used to filter the candidate fluxes on the nonsmooth stencils identified by TENO, ensuring that the filtered nonsmooth flux is within the bounds computed by the MP limiter. In terms of the TVD-type limiters, a second-order (fifth-order) scheme with the Van Albada limiter (the TVD limiter) can be constructed based on the fixed three-cell (five-cell) stencil centered at $i$. Once a candidate stencil is identified as nonsmooth by the TENO weighting strategy, the nonsmooth flux will be replaced by the predefined second-order (fifth-order) flux. Both strategies yield flexible and controllable dissipation for the resulting TENO-M schemes.

Note that, in principle, the TVD-type limiter can also be applied to each identified nonsmooth candidate stencils and the filtering procedure will be similar to that of the MP limiter. The present choice, however, reduces the computational operations since only one limited flux is computed on fixed cells and provides similar performance.}

\subsection{Full high-order construction}
Since the candidate reconstructions are either smooth or filtered by nonlinear limiters, the high-order reconstruction at the cell interface $x_{i+1/2}$ can be achieved by a linear combination
\begin{equation}
\label{eq:reconstuction}
\hat f_{i + 1/2}^K = \sum\limits_{k = 0}^{K - 3} {{w_k}} {\hat f^{\rm{M}}_{k,i + 1/2}}, {\text{ }}k = 0,\cdots,K - 3 ,
\end{equation}
where $w_k = d_k$, which denotes the optimal linear weights. A weight renormalization procedure as with the standard TENO paradigm, Eq.~(\ref{eq:weight_renormalize}), is not needed. In practice, optimal weights $d_k$ can either be optimized to approach the maximum accuracy order or to achieve better spectral properties with low dissipation and low dispersion.

\subsection{Spectral property}
\begin{figure}[htbp]%
\centering
  \includegraphics[width=\textwidth]{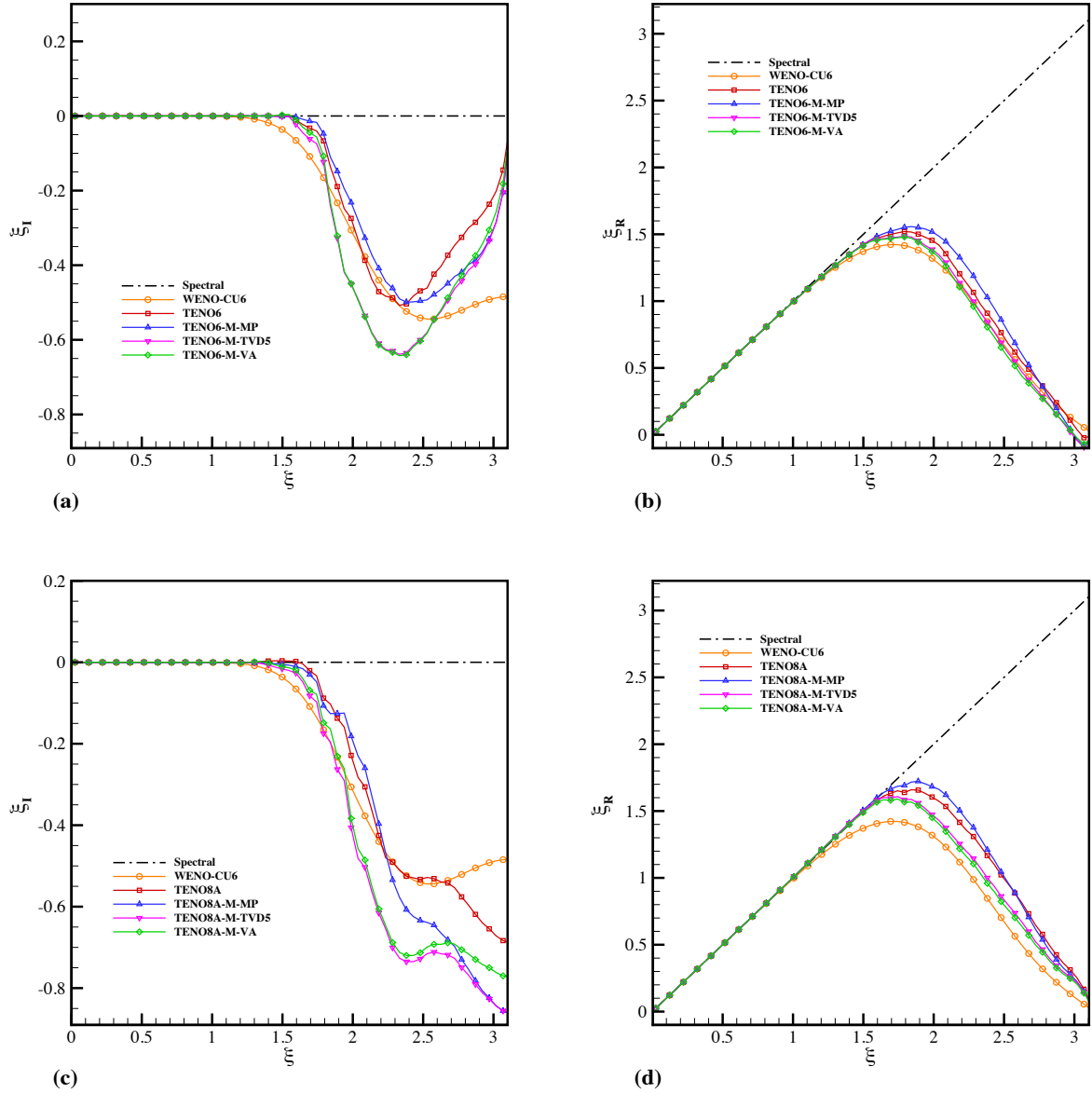}
\caption{Spectral properties. Top: dissipation (a) and dispersion (b) properties of the six-point TENO and TENO-M schemes; bottom: dissipation (c) and dispersion (d) properties of the eight-point TENO and TENO-M schemes.}
\label{Fig:spectral}
\end{figure}

The spectral properties of the nonlinear shock-capturing schemes, i.e. dissipation and dispersion property, can be analyzed by the approximate dispersion relation (ADR) \cite{pirozzoli2006spectral}. As shown in Fig.~\ref{Fig:spectral}, the numerical dissipation and dispersion of the newly proposed six- and eight-point TENO-M schemes are evaluated by computing the imaginary and real parts of the modified wavenumber. The dissipation and dispersion spectra of the standard WENO-CU6 \cite{Hu2010}, TENO6 \cite{fu2016family}, and TENO8A \cite{fu2018targeted} schemes are also plotted for comparison.

It is observed that the good spectral resolutions of TENO6 and TENO8A are generally maintained by both TENO6-M and TENO8A-M schemes, which show remarkable improvements over the WENO-CU6 scheme. Both TENO and TENO-M schemes recover the dispersion and dissipation property of the background optimal linear schemes up to the wavenumber of about 1.5. The adjustment of nonlinear dissipation of TENO-M schemes are restricted to high-wavenumber regions. Owing to the flexibility of the present framework to {choose} nonlinear limiters, the dissipation property of TENO-M schemes in the high-wavenumber regions can be further manipulated by deploying different nonlinear limiters. For the limiters in the present paper, the low- and high-order TVD-type limiters provide more dissipation than the standard TENO schemes. As for the MP limiter, the numerical dissipation of TENO-M-MP schemes can be less than that of the counterpart TENO schemes in certain regions. These results demonstrate that with different choices of nonlinear limiters, the newly proposed framework conveys a broad range of nonlinear-dissipation variations in high wavenumber regions. By the proposed extension, TENO-M schemes can be tailored to handle different kinds of physical computations, such as shocks and turbulence, by a suitable choice of nonlinear limiter.

\section{Numerical validations}

In this section, a set of benchmark cases involving strong discontinuities and broadband flow scales is simulated. The proposed TENO-M schemes are extended to multi-dimensional problems in a dimension-by-dimension manner. For systems of hyperbolic conservation laws, the characteristic decomposition method based on the Roe average \cite{roe1981approximate} is employed. The Rusanov scheme \cite{Rusanov1961} is adopted for flux splitting. The third-order strong stability-preserving (SSP) Runge-Kutta method \cite{gottlieb2001strong} with a typical CFL number of 0.4 is adopted for the time advancement. Numerical results from WENO-CU6 \cite{Hu2010}, TENO6 \cite{fu2016family} and TENO8A \cite{fu2018targeted} are compared. The parameters in TENO6 and TENO8A schemes are identical to those in \cite{fu2016family} and \cite{fu2018targeted}.

\subsection{Advection problems}
\subsubsection{Accuracy test}

We first consider the one-dimensional Gaussian pulse advection problem \cite{yamaleev2009systematic}. The linear advection equation
\begin{equation}
\label{eq:linear_advect}
\frac{\partial u}{\partial t} + \frac{\partial u}{\partial x} = 0 ,
\end{equation}
with initial condition
\begin{equation}
\label{eq:linear_advect_init}
{u}(x,0) = {e^{ - 300(x - {x_c})^2}}, \text{ }  x_c = 0.5 ,
\end{equation}
is solved in a computational domain $ 0 \leq x \leq 1$ and the final time is $t = 1$. Periodic boundary conditions are imposed at $ x = 0$ and $x = 1$.

\begin{figure}[htbp]%
\centering
  \includegraphics[width=0.95\textwidth]{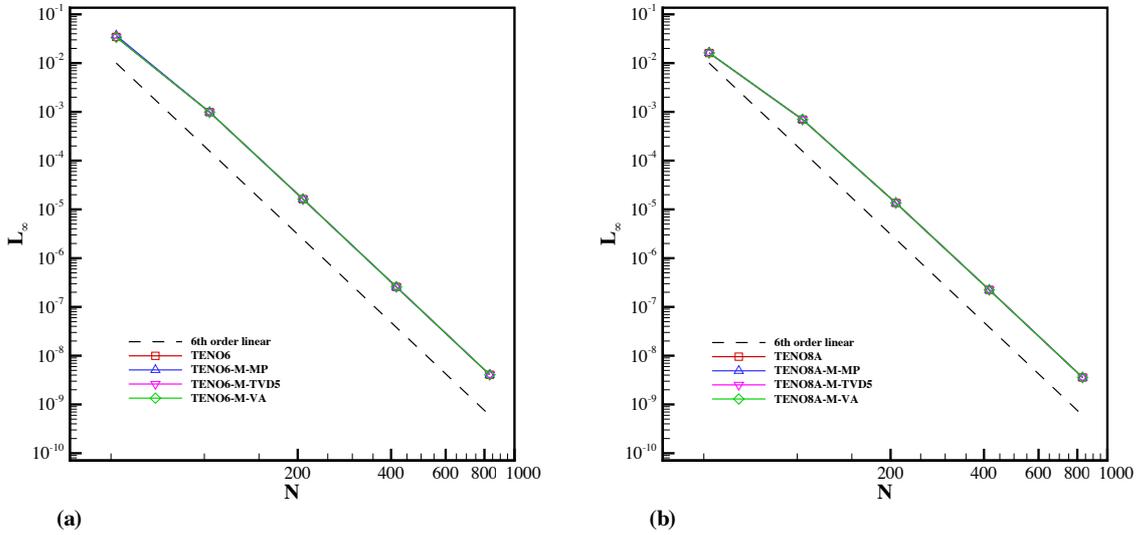}
\caption{Convergence of the $L_\infty$ error from the TENO6-M (a) and TENO8A-M (b) schemes. $N$ denotes the grid number.}
\label{Fig:Convergence}
\end{figure}

As shown in Fig.~\ref{Fig:Convergence}, the desired accuracy order is achieved for both TENO6-M and TENO8A-M schemes. As all solutions is identified as smooth, optimal accuracy order is achieved.

\subsubsection{Multi-wave advection}
This case is taken from Jiang and Shu \cite{Jiang1996}. We solve the advection equation
\begin{equation}
\label{eq:linear_advect}
\frac{\partial u}{\partial t} + \frac{\partial u}{\partial x} = 0 ,
\end{equation}
with the initial condition
\begin{equation}
\label{eq:threepart}
{u}(x,0) = \left\{ {\begin{array}{*{20}{c}}
{\frac{1}{6}[G(x-1, \beta, z-\theta)+ G(x-1, \beta, z+\theta) + 4G(x-1, \beta, z)],} &{\text{if }0.2 \le x < 0.4} ,\\
{1,}&{\text{if }0.6\le x \le 0.8} ,  \\
{1-|10(x-1.1)|,}&{\text{if }1.0\le x \le 1.2} ,  \\
{\frac{1}{6}[F(x-1, \alpha, a-\theta)+ F(x-1, \alpha, a+\theta) + 4F(x-1, \alpha, a)],} &{\text{if } 1.4 \le x <
1.6} ,\\
{0,} &{\text{otherwise,}}
\end{array}} \right.
\end{equation}
where
\begin{equation}
\label{eq:intialsin}
G(x, \beta, z) = e^{-\beta(x-z)^2} , \\
F(x, \alpha, a) = \sqrt{\max(1-\alpha^2(x-a)^2, 0)} .
\end{equation}
The parameters in Eqs.~(\ref{eq:threepart}--\ref{eq:intialsin}) are given as
\begin{equation}
\label{eq:constant}
a = 0.5, \text{ } z = -0.7, \text{ } \theta =0.005, \text{ } \alpha =10, \text{ } \beta = \frac{\log2}{36\theta^2}.
\end{equation}
The initial distribution consists of a Gaussian pulse, a square wave, a sharply peaked triangle and a half-ellipse arranged from left to right in the computational domain $x\in[0,2]$. The computation is performed with $N = 200$ uniformly distributed mesh cells and the final simulation time is $t = 2$.

\begin{figure}[htbp]%
\centering
  \includegraphics[width=0.95\textwidth]{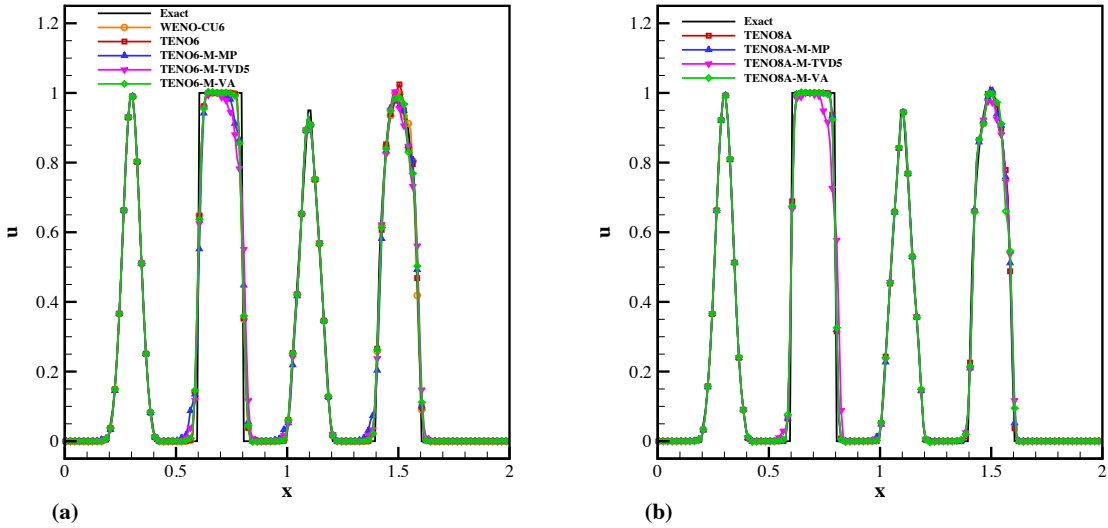}
\caption{Numerical results of advection of multi-wave with TENO6-M (a) and TENO8A-M (b) schemes. Discretization  is on 200 uniformly distributed grid points. }
\label{Fig:full_multiwave}
\end{figure}

As shown in Fig.~\ref{Fig:full_multiwave}, for the advection of the square wave, all schemes are free from numerical oscillations. Concerning the advection of the half-ellipse wave, TENO-M schemes eliminate the overshoots generated by the standard TENO and WENO-CU6 schemes. Qualitatively, TENO-M-VA schemes exhibit best symmetry preservation.

\subsection{Gas dynamics}

\subsubsection{Shock-tube problem}
The Lax problem \cite{Lax1954} and the Sod problem \cite{Sod1978} are considered.
The initial condition for Lax problem \cite{Lax1954} is
\begin{equation}
\label{eq:lax}
(\rho ,u,p) = \left\{ {\begin{array}{*{20}{c}}
{(0.445,0.698,3.528),}&{\text{if }0 \le x < 0.5} ,\\
{(0.5,0,0.5710),}&{\text{if }0.5 \le x \le 1} ,
\end{array}} \right.
\end{equation}
and the final simulation time is $t = 0.14$.

The initial condition for Sod problem \cite{Sod1978} is
\begin{equation}
\label{eq:sod}
(\rho ,u,p) = \left\{ {\begin{array}{*{20}{c}}
{(1,0,1),}&{\text{if }0 \le x < 0.5} ,\\
{(0.125,0,0.1),}&{\text{if }0.5 \le x \le 1} ,
\end{array}} \right.
\end{equation}
and the final simulation time is $t = 0.2$. The computations are done on 100 uniformly distributed grid points.

\begin{figure}[htbp]%
\centering
  \includegraphics[width=0.95\textwidth]{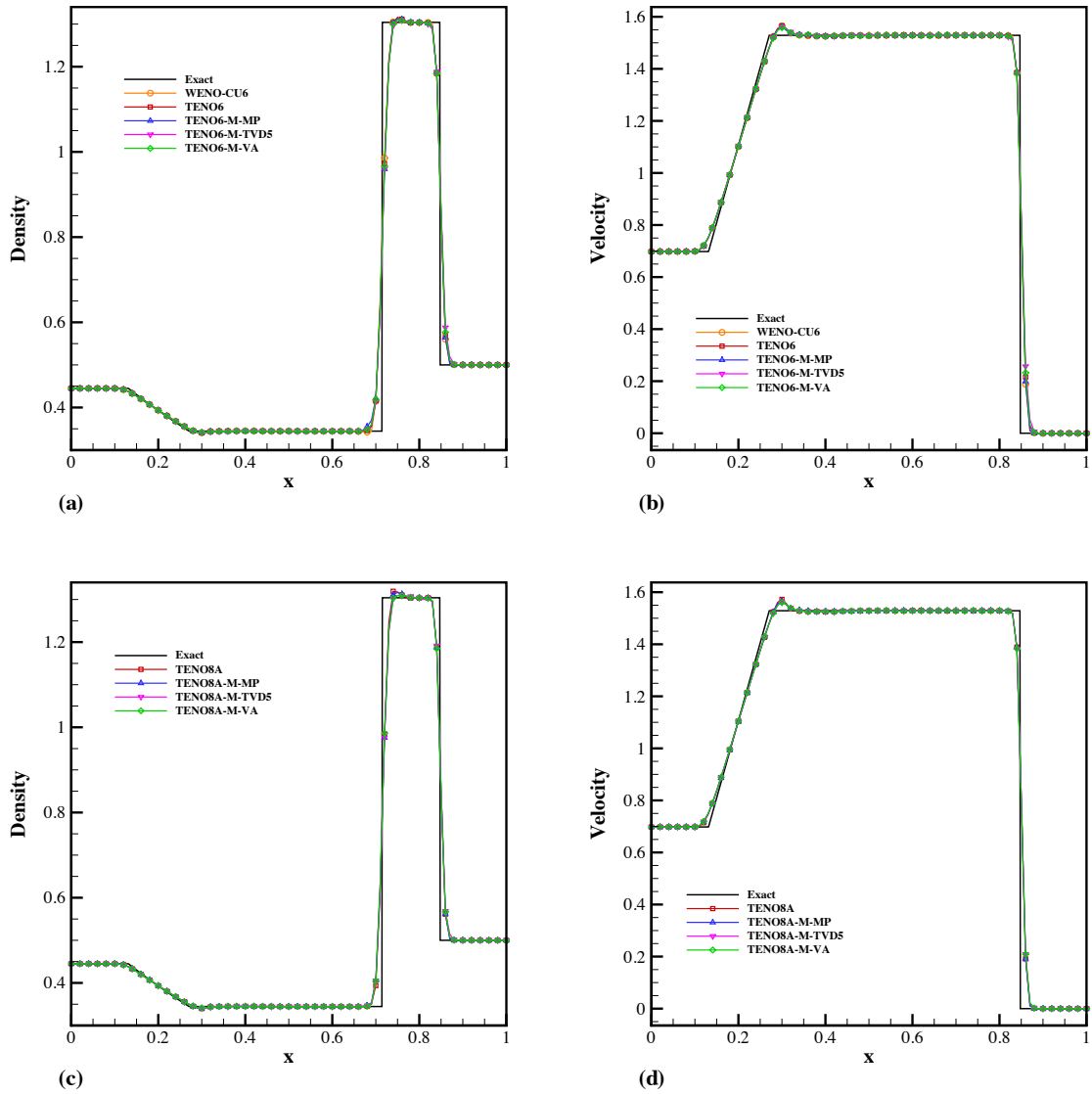}
\caption{Shock-tube problem: the Lax problem. Top: density profile (a) and velocity profile (b) of WENO-CU6, TENO6 and TENO6-M schemes; bottom: density profile (c) and velocity profile (d) of the TENO8A and TENO8A-M schemes. Discretization is on 100 uniformly distributed grid points and the final simulation time is $t = 0.14$.}
\label{Fig:lax}
\end{figure}

\begin{figure}[htbp]%
\centering
  \includegraphics[width=0.95\textwidth]{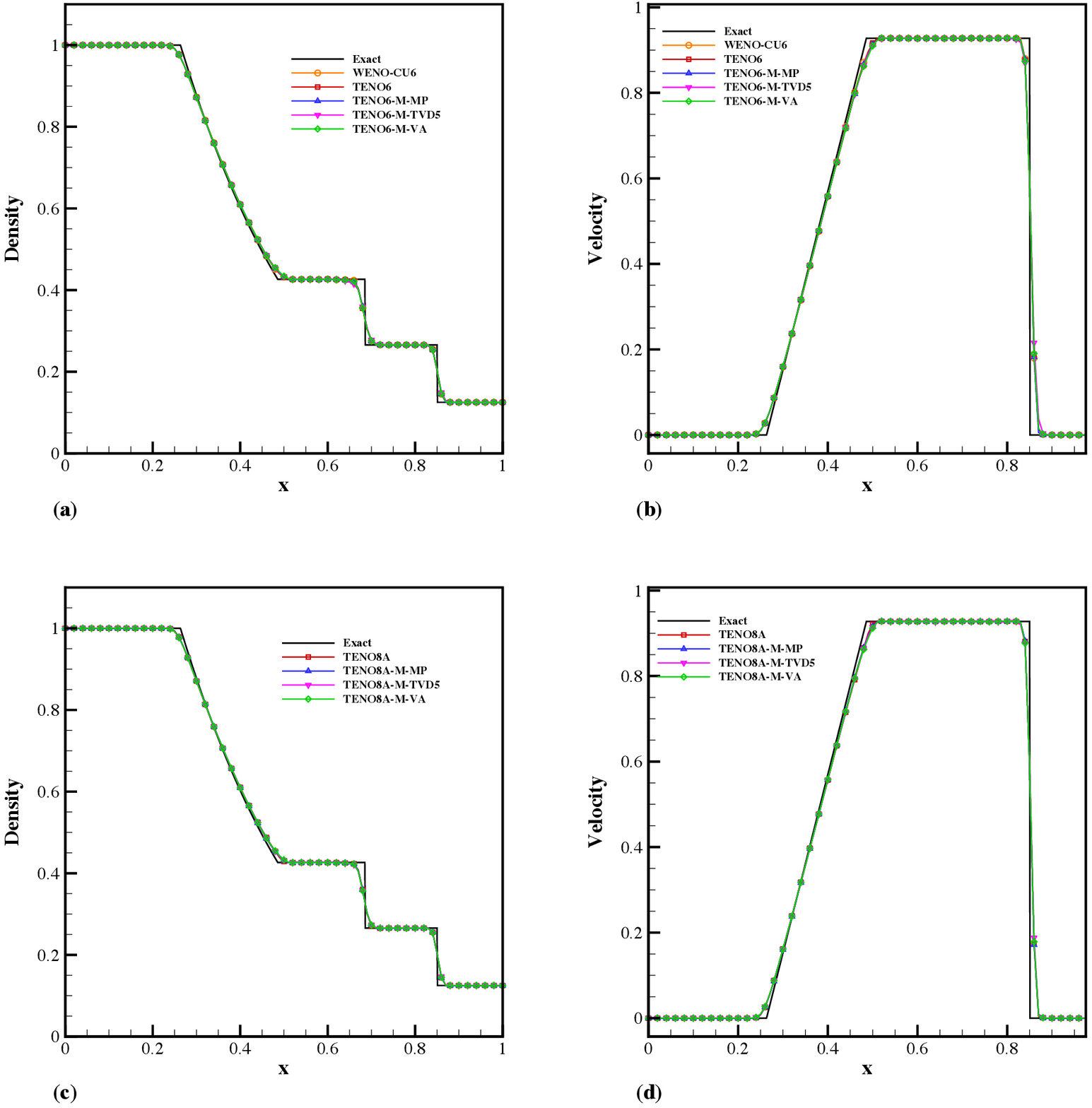}
\caption{Shock-tube problem: the Sod problem. Top: density profile (a) and velocity profile (b) of WENO-CU6, TENO6 and TENO6-M schemes; bottom: density profile (c) and velocity profile (d) of TENO8A and TENO8A-M schemes. Discretization is on 100 uniformly distributed grid points and the final simulation time is $t = 0.2$. }
\label{Fig:sod}
\end{figure}

As shown in Fig.~\ref{Fig:lax} and Fig.~\ref{Fig:sod}, both the proposed TENO6-M and TENO8A-M schemes recover the shock-capturing capability of standard TENO schemes {with the nonsmooth stencils being filtered} by nonlinear limiters.
\subsubsection{Shock density-wave interaction problem}
This case is proposed by Shu and Osher \cite{Shu1989}. A one-dimensional Mach-3 shock wave interacts with a perturbed density field generating both small-scale structures and discontinuities. The initial condition is
\begin{equation}
\label{eq:shuosher}
(\rho ,u,p) = \left\{ {\begin{array}{*{20}{c}}
{(3.857,2.629,10.333),}&{\text{if } 0 \le x < 1} ,\\
{(1 + 0.2\sin (5(x-5)),0,1),}&{\text{if } 1 \le x \le 10} .
\end{array}} \right.
\end{equation}
The computational domain is [0,10] with $N = 200$ uniformly distributed mesh cells and the final evolution time is $t = 1.8$. The exact solution for reference is obtained by the fifth-order WENO-JS scheme with $N = 2000$.

\begin{figure}[htbp]%
\centering
  \includegraphics[width=0.95\textwidth]{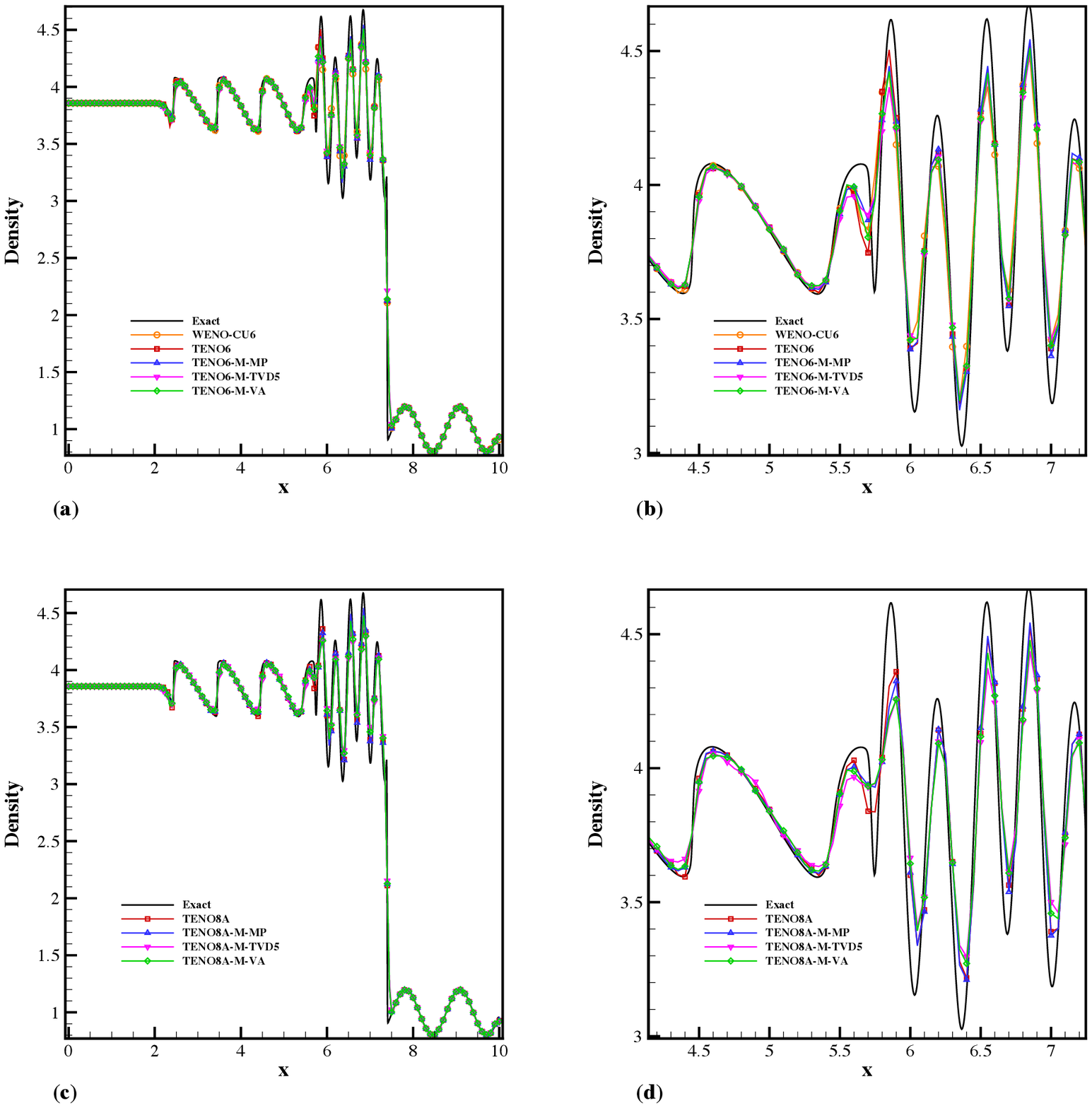}
\caption{Shock density-wave interaction problem: solutions from the WENO-CU6, TENO and TENO-M schemes. Top: density distribution (a) and a zoom-in view of the density distribution (b) of WENO-CU6,  six-point TENO6 and TENO6-M schemes; bottom: density distribution (c) and a zoom-in view of the density distribution (d) of eight-point TENO8A and TENO8A-M schemes. Discretization is on 200 uniformly distributed grid points and the final simulation time is $t = 1.8$. }
 \label{Fig:shu-osher}
\end{figure}

Fig.~\ref{Fig:shu-osher} shows the computed density distributions from both six- and eight-point TENO-M schemes. It is observed that the proposed TENO6-M and TENO8A-M schemes resolve the high-wavenumber fluctuations with very low numerical dissipation while capturing the shocklets as sharply as the WENO-CU6 and standard TENO schemes. Among TENO-M schemes, the TENO-M-TVD5 and TENO-M-VA schemes provide slightly more dissipation than the standard TENO schemes in the under-resolved regions while the TENO-M-MP schemes are less dissipative with better wave-resolution property. By adapting the nonlinear limiters, the numerical dissipation in nonsmooth regions is tuned accordingly without sacrificing the shock-capturing capability and the high-order accuracy with this new TENO framework.

\subsubsection{Interacting blast waves}
The two-blast-waves interaction problem taken from \cite{Woodward1984} is considered. The initial condition is
\begin{equation}
\label{eq:blastwaves}
(\rho ,u,p) = \left\{ {\begin{array}{*{20}{c}}
{(1,0,1000),}&{\text{if }0 \le x < 0.1},\\
{(1,0,0.01),}&{\text{if }0.1 \le x < 0.9},\\
{(1,0,100),}&{\text{if }0.9 \le x \le 1} .
\end{array}} \right.
\end{equation}
The simulation is performed on a uniform mesh with $N = 400$ and the final simulation time is $t = 0.038$. The exact solution for reference is computed by the fifth-order WENO-JS scheme on a uniform mesh with $N = 2500$. For this case the Roe scheme with entropy-fix is employed as numerical flux function.

\begin{figure}[htbp]%
\centering
  \includegraphics[width=0.95\textwidth]{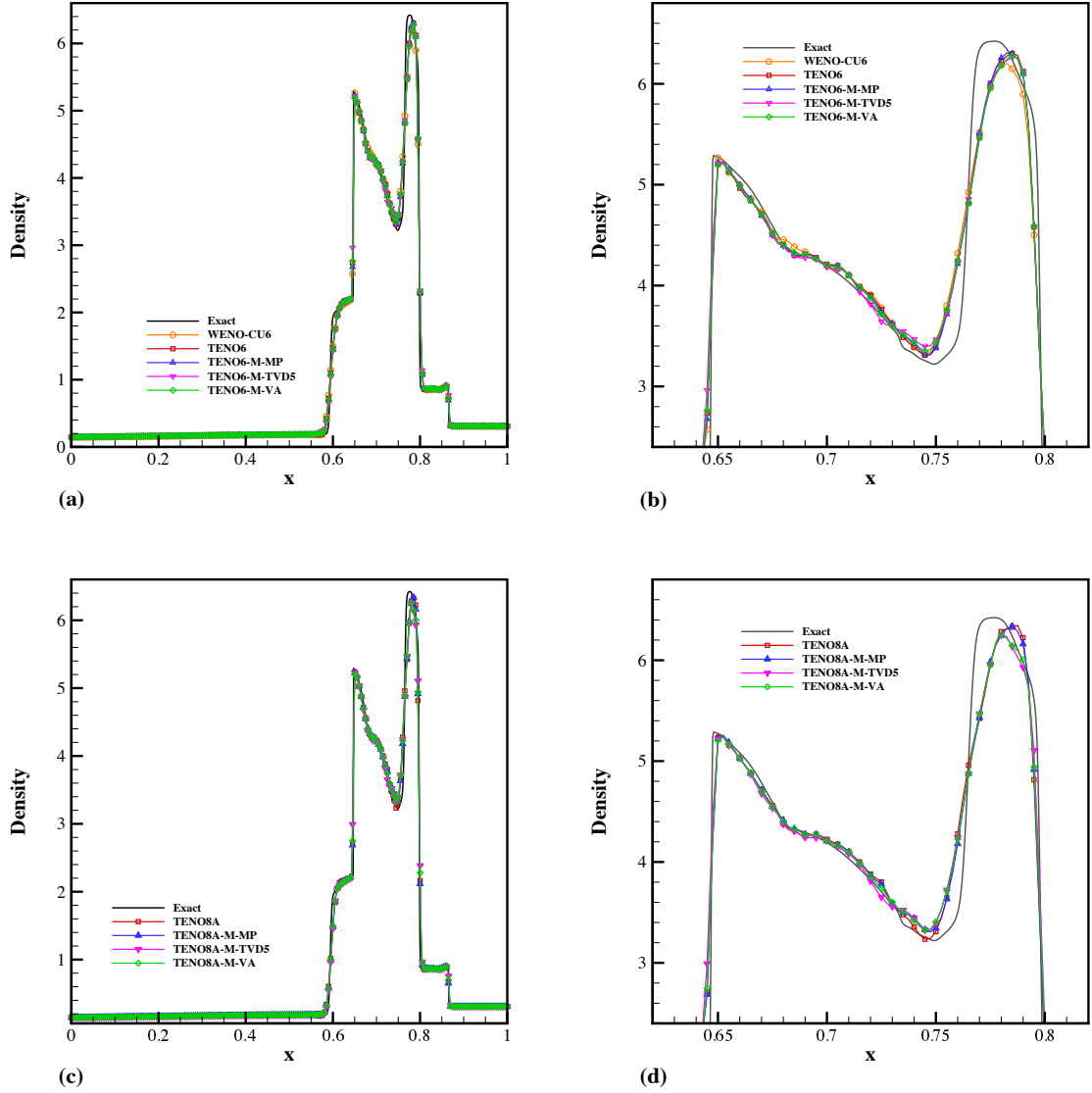}
\caption{Interacting blast waves problem: solutions from the WENO-CU6, TENO and TENO-M schemes. Top: density distribution (a) and a zoom-in view of the density distribution (b) of WENO-CU6, six-point TENO and TENO-M schemes; bottom: density distribution (c) and a zoom-in view of the density distribution (d) of the eight-point TENO and TENO-M schemes. Discretization is on 400 uniformly distributed grid points and the final simulation time is $t = 0.038$.}
\label{Fig:blastwaveproblem}
\end{figure}

As shown in Fig.~\ref{Fig:blastwaveproblem}, all schemes generate good results and the TENO6-M and TENO8A-M schemes perform better than WENO-CU6 in resolving the density peak at $x = 0.78$. In addition, TENO-M-MP schemes provide similar dissipation compared to the standard TENO schemes while TENO-M-TVD5 and TENO-M-VA schemes provide more dissipation, which confirms that a broad range of nonlinear-dissipation variations can be achieved through adopting different nonlinear limiters.
\subsubsection{Rayleigh-Taylor instability}

The inviscid Rayleigh-Taylor instability case proposed by Xu and Shu \cite{Xu2005} is considered here. The initial condition is
\begin{equation}
(\rho ,u,v,p) = \left\{ {\begin{array}{*{20}{c}}
{(2,0, - 0.025c\cos (8\pi x),1 + 2y),}&{\text{if } 0 \le y < 1/2},\\
{(1,0, - 0.025c\cos (8\pi x),y + 3/2),}&{\text{if } 1/2 \le y \le 1},
\end{array}} \right.
\end{equation}
where the sound speed $c = \sqrt {\gamma \frac{p}{\rho }}$ with $\gamma  = \frac{5}{3}$. The computational domain is $[0,0.25] \times [0,1]$. Reflective boundary conditions are imposed at the left and right side of the domain. Constant primitive variables $(\rho ,u,v,p) = (2,0,0,1)$ and $(\rho ,u,v,p) = (1,0,0,2.5)$ are set for the bottom and top boundaries, respectively.

\begin{figure}[htbp]%
\centering
\includegraphics[width=\textwidth]{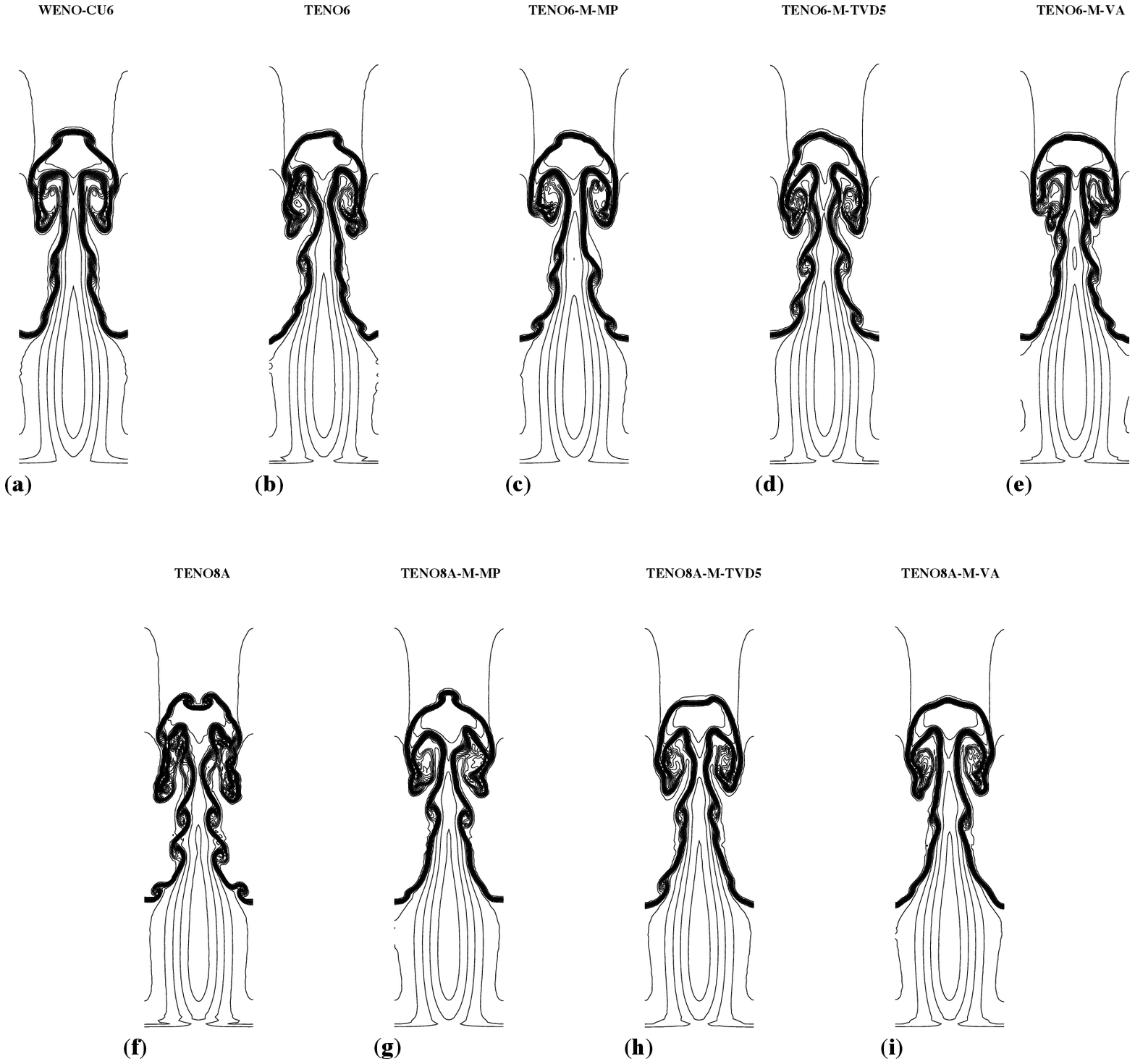}
  \caption{Rayleigh-Taylor instability: density contours from the WENO-CU6 (a), TENO6 (b), TENO6-M (c)(d)(e), TENO8A (f) and TENO8A-M (g)(h)(i) schemes at simulation time $t = 1.95$. Resolution is $64 \times 256$. This figure is drawn with 25 density contours between 0.9 and 2.2.}
 \label{Fig:RTI}
\end{figure}

In Fig.~\ref{Fig:RTI}, the computed density contours with the proposed TENO6-M and TENO8A-M schemes at resolution of $64 \times 256$ are shown.
The results show that the newly proposed TENO6-M and TENO8A-M schemes resolve finer small-scale structures than WENO-CU6 scheme. For all newly proposed schemes, the nonsmooth contact interface regions are captured sharply. In terms of vortex structures, TENO6-M-MP and TENO6-M-TVD5 schemes generate finer structures compared with TENO6-M-VA. While TENO6-M-MP exhibits breaking of flow symmetry, better symmetry preservation is visible for TENO6-M-VA. Similar observations can be seen for TENO8A and TENO8A-M.

\subsubsection{Double Mach reflection of a strong shock}
This two-dimensional case is taken from Woodward and Colella \cite{Woodward1984}. The initial condition is
\begin{equation}
\label{eq:DMR}
(\rho ,u,v,p) = \left\{ {\begin{array}{*{20}{c}}
{(1.4,0,0,1) ,}&{\text{if } y < 1.732(x - 0.1667) ,}\\
{(8,7.145, - 4.125,116.8333) , }&{\text{otherwise}.}
\end{array}} \right.
\end{equation}
The computational domain is $[0,4] \times [0,1]$ and the final simulation time is $t = 0.2$. The boundary condition setup is the same as \cite{Woodward1984}.

\begin{figure}[htbp]
\centering
\includegraphics[width=0.5\textwidth]{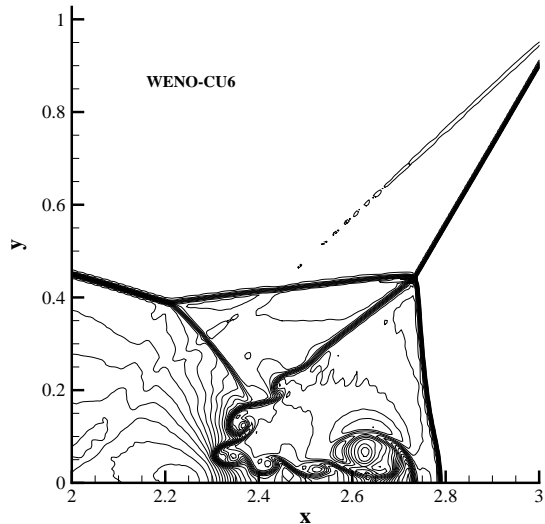}
  \caption{Double Mach reflection of a strong shock: density contours from WENO-CU6 for reference at simulation time $t = 0.2$. Resolution is $800 \times 200$. This figure is drawn with 42 density contours between 3.27335 and 20.1335.}
 \label{Fig:dmr1024wenocu}
\end{figure}

\begin{figure}[htbp]%
\centering
\includegraphics[width=0.9\textwidth]{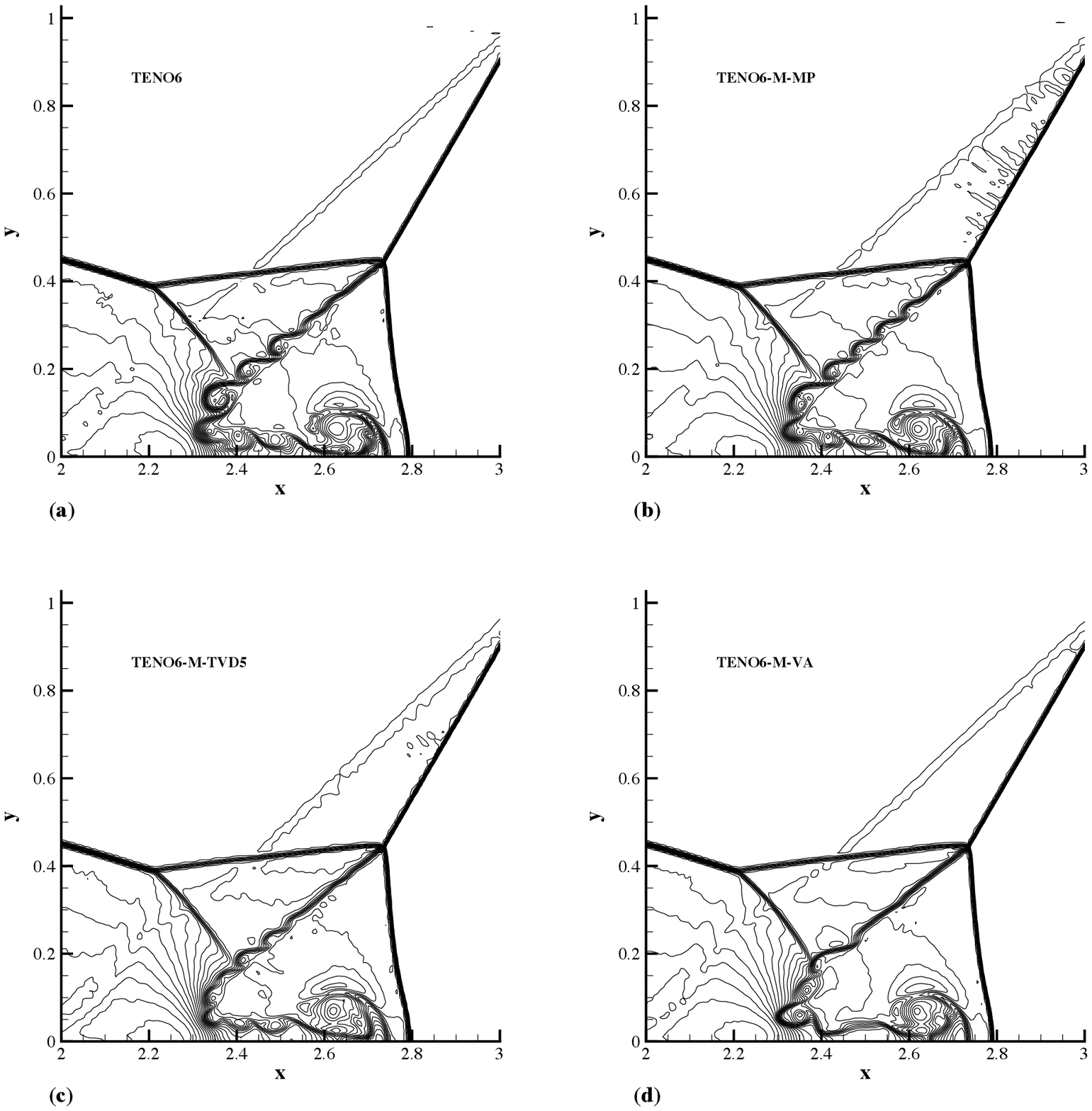}
  \caption{Double Mach reflection of a strong shock: density contours from TENO6 (a) and TENO6-M schemes (b)(c)(d) at simulation time $t = 0.2$. Resolution is $800 \times 200$. This figure is drawn with 42 density contours between 3.27335 and 20.1335.}
 \label{Fig:dmr1024DMR6th}
\end{figure}

\begin{figure}[htbp]%
\centering
\includegraphics[width=0.9\textwidth]{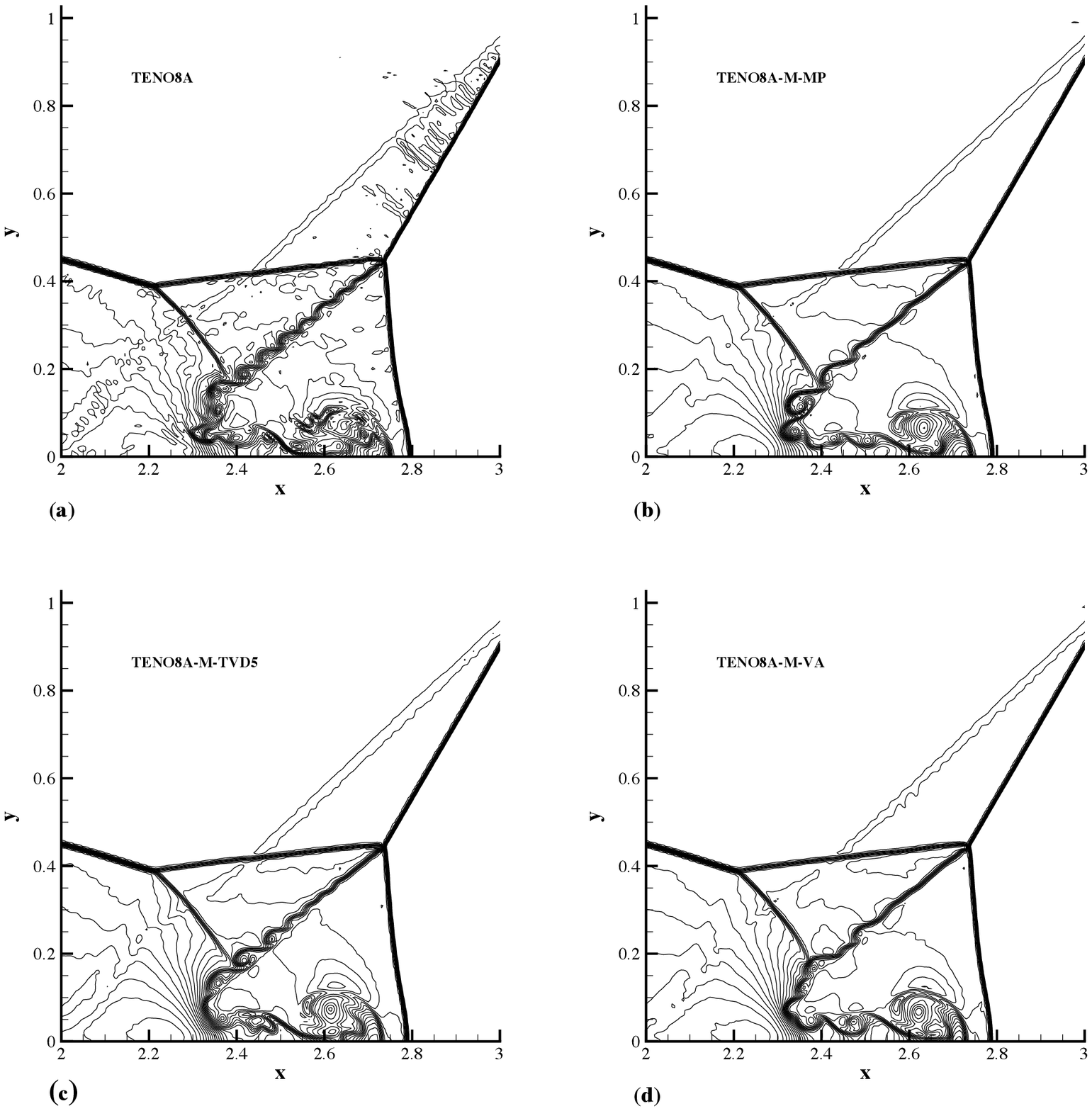}
  \caption{Double Mach reflection of a strong shock: density contours from TENO8A (a) and TENO8A-M schemes (b)(c)(d) at simulation time $t = 0.2$. Resolution is $800 \times 200$. This figure is drawn with 42 density contours between 3.27335 and 20.1335.}
 \label{Fig:dmr1024DMR8Ath}
\end{figure}

As shown in Fig.~\ref{Fig:dmr1024wenocu}, Fig.~\ref{Fig:dmr1024DMR6th} and Fig.~\ref{Fig:dmr1024DMR8Ath}, TENO and TENO-M schemes perform better than WENO-CU6 in resolving the small-scale structures. The shock-wave patterns are captured sharply without spurious numerical oscillations. Compared with the results of the TENO6-M-MP scheme, that of TENO6-M-TVD5 and TENO6-M-VA are slightly more dissipative.

\begin{figure}[htbp]%
\centering
\includegraphics[width=0.9\textwidth]{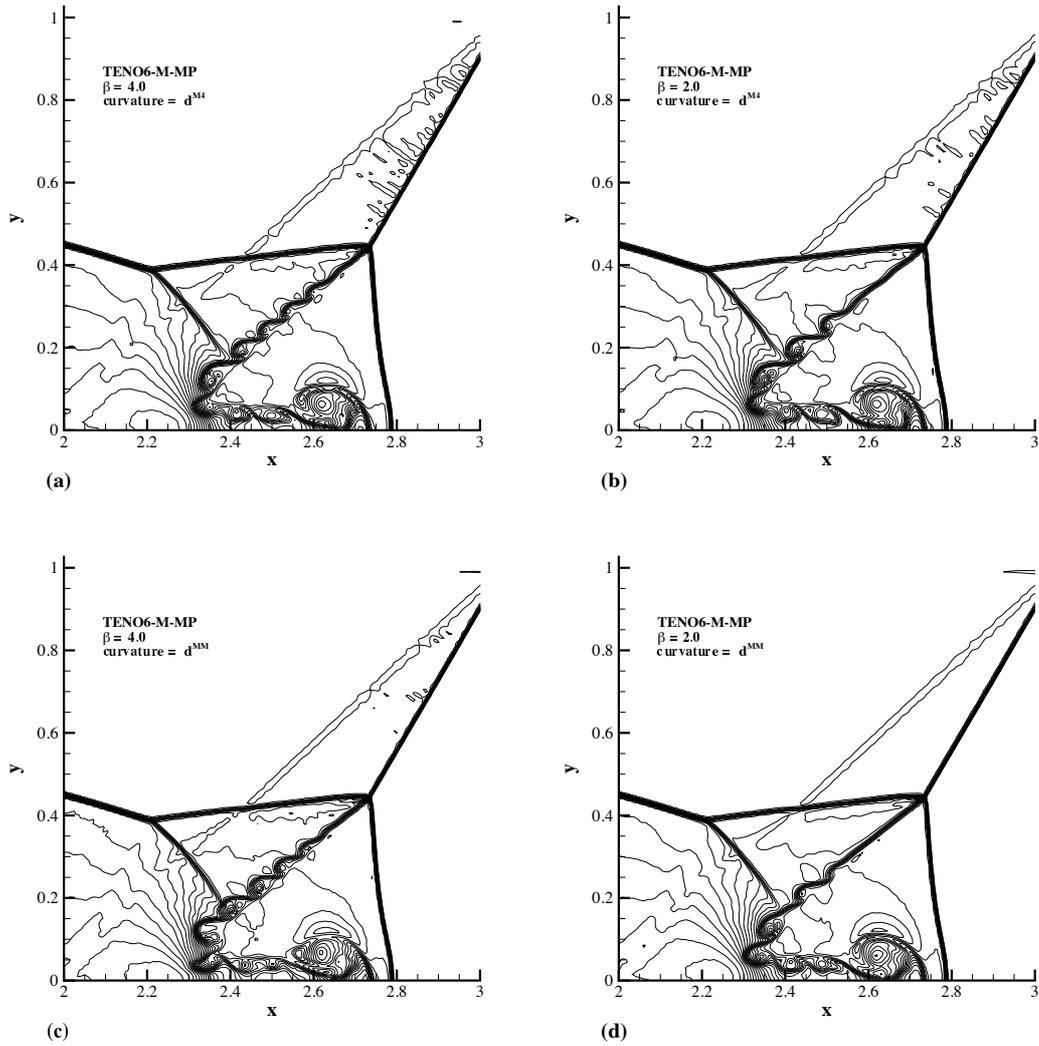}
  \caption{Double Mach reflection of a strong shock: parameter study of monotonicity-preserving limiter.  Density contours from TENO6-M-MP with $d_{M4}$ curvature measurement and $\beta = 4.0$ (a), TENO6-M-MP with $d_{M4}$ curvature measurement and $\beta = 2.0$ (b), TENO6-M-MP with $d_{MM}$ curvature measurement and $\beta = 4.0$ (c) and TENO6-M-MP with $d_{MM}$ curvature measurement and $\beta = 2.0$ (d) at simulation time $t = 0.2$. Resolution is $800 \times 200$. This figure is drawn with 42 density contours between 3.27335 and 20.1335.}
 \label{Fig:parameterstudy}
\end{figure}

Moreover, to address the flexibility of the unified framework, {it is worth noting that the effect of the MP limiter can be adjusted through tuning of the curvature measurement $d_{i + 1/2}^{\text{M4}}$ and $\beta$. }
In the following, we test the dissipation property of MP limiter with a smaller $\beta = 2$ and a less restrictive curvature measurement $d_{i + 1/2}^{\text{MM}}$ at the cell interface $i+ 1/2$ which is defined as
\begin{equation}
\label{eq:maximal_curvature1}
{d_{i + 1/2}^{\text{MM}}} = {\rm{minmod}} ({d_i},{d_{i + 1}}).
\end{equation}
As shown in Fig.~\ref{Fig:parameterstudy}, the dissipation property is more sensitive to $\beta$, which determines the amount of freedom allowing for large curvature. And Fig.~\ref{Fig:parameterstudy} shows that, a more restrictive curvature measurement results in a less dissipative MP limiter.

\begin{figure}[htbp]
\centering
\includegraphics[width=0.45\textwidth]{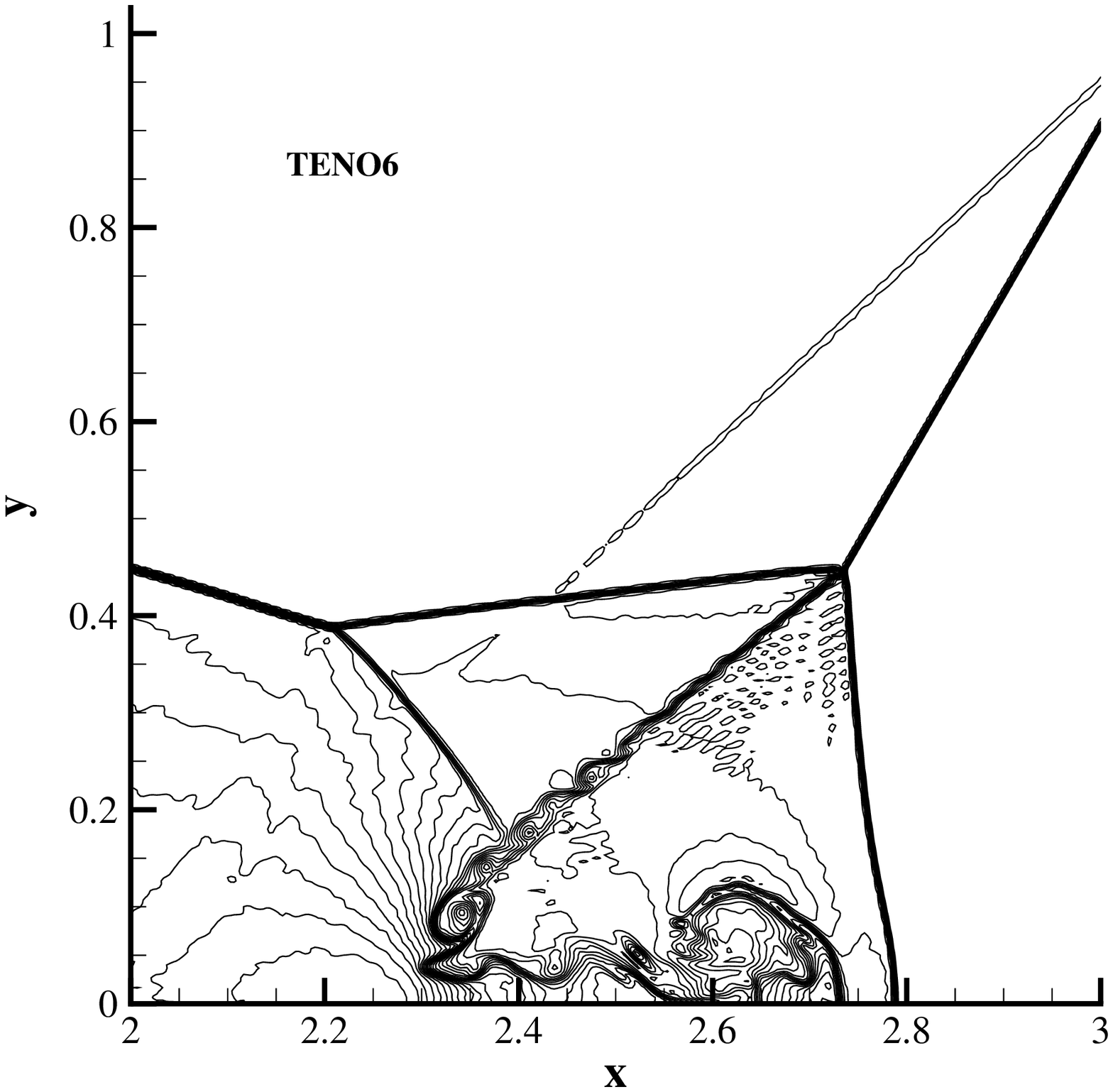}
  \caption{Double Mach reflection of a strong shock: density contours from TENO6 at simulation time $t = 0.2$. Resolution of $1200 \times 300$. This figure is drawn with 42 density contours between 3.27335 and 20.1335.}
 \label{Fig:DMRteno6RoeM}
\end{figure}

\begin{figure}[htbp]
\centering
\includegraphics[width=0.9\textwidth]{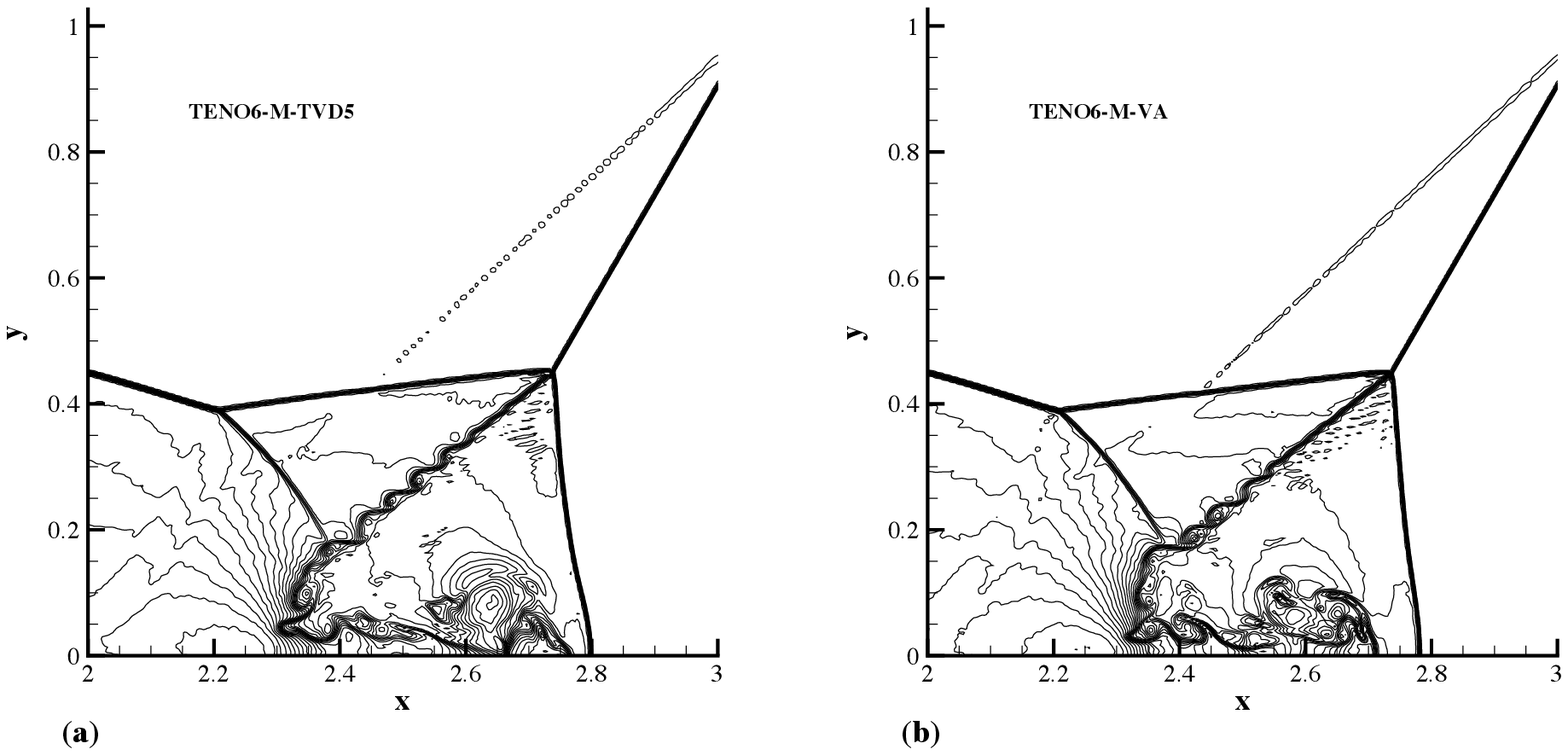}
  \caption{Double Mach reflection of a strong shock: density contours from TENO6-M schemes: TENO6-M-TVD5 (a) and TENO6-M-VA (b) at simulation time $t = 0.2$. TENO6-M-MP fails to pass. Resolution of $1200 \times 300$. This figure is drawn with 42 density contours between 3.27335 and 20.1335.}
 \label{Fig:DMRteno6MRoeM}
\end{figure}

We further perform simulations with a low-dissipation flux-splitting method and higher resolution. In Fig.~\ref{Fig:DMRteno6RoeM} and Fig.~\ref{Fig:DMRteno6MRoeM}, simulations with the resolution of $N = 1200 \times 300$ are performed on uniform meshes with TENO6 and TENO6-M schemes. Following Fleischmann et.al \cite{fleischmann2019low}, to cure the grid-aligned shock instability in low-dissipation computations, the newly proposed Roe-M flux-splitting method is deployed to perform the simulations. With additional nonlinear limiters to filter the nonsmooth candidate stencils, TENO6-M-TVD5 and TENO6-M-VA can pass the case without difficulties. As shown in Fig.~\ref{Fig:DMRteno6MRoeM}, TENO6-M schemes capture the shock-wave pattern sharply. In terms of the vortex structures, TENO6-M-TVD5 and TENO6-M-VA schemes generate these small scale structures with reduced numerical noise.



\subsection{Turbulent flows}
\begin{figure}[htbp]
\centering
\subfigure{\includegraphics[width=0.95\textwidth]{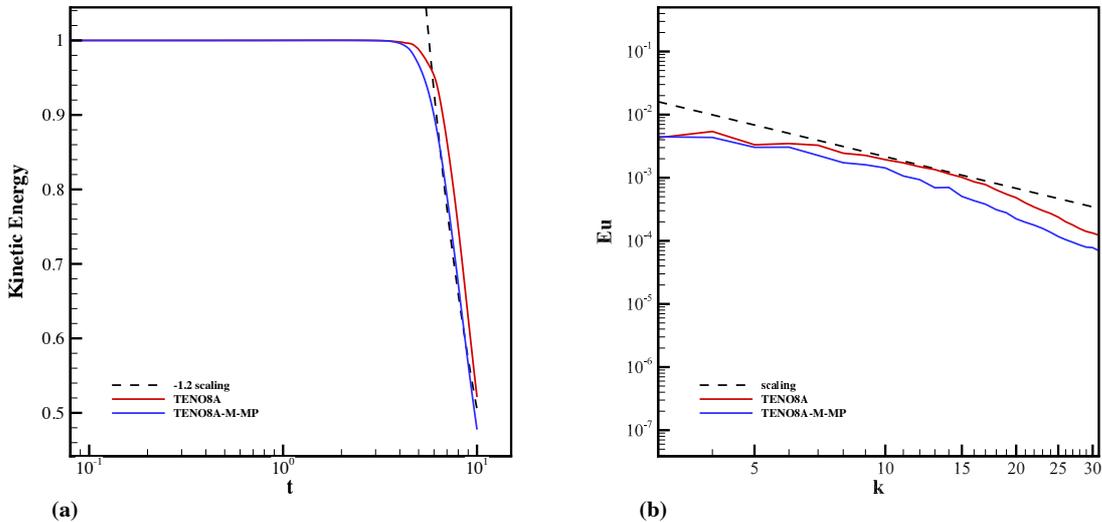}}
\caption{Inviscid Taylor-Green vortex problem: solutions from the TENO8A and TENO8A-M schemes. The evolution of normalized total kinetic energy (a) and the energy spectrum within resolved inertial subrange at t = 10 compared to the Kolmogorov scaling $E(k) \propto k^{-5/3}$ (b). Discretization is on $64^3$ uniformly distributed grid points.}
\label{Fig:inviscidTGV}
\end{figure}

{
\subsubsection{Inviscid incompressible isotropic Taylor-Green vortex (TGV)}

In order to assess the proposed TENO-M scheme for implicit large eddy simulations (ILES), we consider the typical incompressible isotropic Taylor-Green vortex case. The initial condition follows
\begin{equation}
\label{eq:TGV}
{\begin{array}{*{20}{c}}
{u(x,y,z,0)}=&{\text{sin(x)cos(y)cos(z) }},\\
{v(x,y,z,0)}=&{\text{-cos(x)sin(y)cos(z)}},\\
{w(x,y,z,0)}=&{0},\\
{\rho(x,y,z,0)}=&{\text{1.0} },\\
{p(x,y,z,0)}=&{100+\frac{1}{16}[\text{(cos(2x)+cos(2y))(2+cos(2z))}-2]}.\\
\end{array}}
\end{equation}
The computational domain is $[0,2\pi] \times [0,2\pi] \times [0,2\pi]$ with periodic boundary conditions applied on the six boundaries. The mesh resolution is $64^3$. The evolution of normalized total kinetic energy is shown in Fig.~\ref{Fig:inviscidTGV}(a). The kinetic energy decays as $t^{-1.2}$, which agrees well with \cite{lesieur20003d}. For very large Reynolds number incompressible TGV flows, three-dimensional statistically isotropic turbulence develops, and the kinetic-energy spectra $E(k) \propto k^{-5/3}$ is observed within the inertial subrange after $t \simeq 9$. As shown in Fig.~\ref{Fig:inviscidTGV}(b), although TENO8A-M is slightly more dissipative than TENO8A (for which the built-in parameters are particularly optimized for this case \cite{fu2018targeted}) at high wavenumbers, it can faithfully reproduce the Kolmogorov scaling up to 2/3 of the cut-off wavenumber at $t = 10$. The present results suggest that TENO8-A-M scheme can be applied to under-resolved simulations without parameter tuning and has potential to function as an ILES model. }

\begin{figure}[htbp]
\centering
\subfigure{\includegraphics[width=0.95\textwidth]{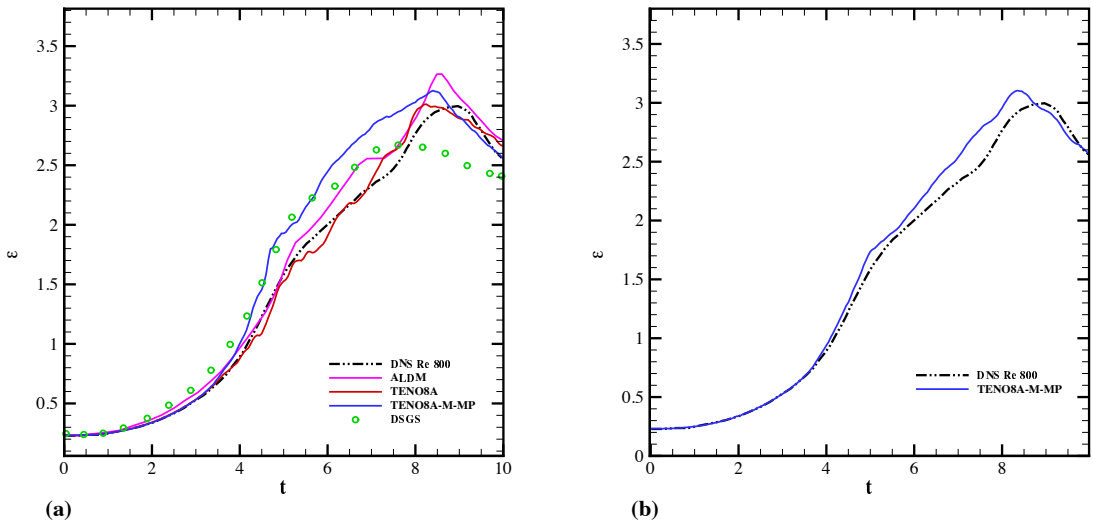}}
\caption{Taylor-Green vortex problem with Reynolds number 800: comparison of energy dissipation rate. Left: comparison on a $64^3$ mesh (a). Right: comparison on a $96^3$ mesh (b).}
\label{Fig:TGV800}
\end{figure}

{
\subsubsection{Viscous incompressible isotropic Taylor-Green vortex}
With finite Reynolds numbers, the favorable numerical method is expected to reproduce the flow transition and the developed turbulence. Viscous incompressible isotropic case with a Reynolds number of 800 is computed with coarse resolutions of $64^3$ and $96^3$. Results from the implicit LES model ALDM \cite{hickel2006adaptive} and the dynamic Smagorinsky model (DSGS) \cite{germano1991dynamic} are also provided for comparisons.
As shown in Fig.~\ref{Fig:TGV800}, TENO8A-M-MP exhibits slightly higher dissipation than ALDM and TENO8A (both are optimized for this case \cite{fu2018targeted}\cite{hickel2006adaptive}) at the time stage from $t = 5$ to $t = 8$ on the coarse mesh of $64^3$. However, the result is better than that of the dynamic Smagorinsky model. With the finer mesh of $96^3$, TENO8A-M-MP scheme shows good agreement with the DNS reference (computed on a mesh of $256^3$ \cite{brachet1983small}) and predicts the peak dissipation rate at the later time stage very well. Therefore, TENO8A-M-MP can be applied as a physically-consistent ILES model without parameter tuning on coarse meshes.
}

{
\subsection{Extreme simulations}
In this part, we consider two challenging cases to demonstrate the advantage of the newly proposed framework. Following \cite{fu2019very}, the positivity-preserving flux limiter \cite{hu2013positivity} is applied to guarantee positivity of pressure and density. With the present MP limiter, a less restrictive curvature measurement $d_{i + 1/2}^{\text{MM}}$ and $\beta = 1$ are adopted to eliminate artificial local extrema generated by spurious oscillations.
\subsubsection{Noh problem}
The Noh problem \cite{noh1987errors} is a challenging case to test shock-capturing algorithms, where a uniform inwards radial inflow generates an infinitely strong shock wave moving outwards at a constant speed. The computational domain is $[0,0.256] \times [0,0.256] \times [0,0.256]$ and the initial condition is described as
 $(\rho ,u,v,w,p) = (1,-x/r,-y/r,-z/r,1\times10^{-6})$. The final simulation time is 0.6 and the mesh resolution is $64^3$.

\begin{figure}[htbp]
\centering
\subfigure{\includegraphics[width=0.95\textwidth]{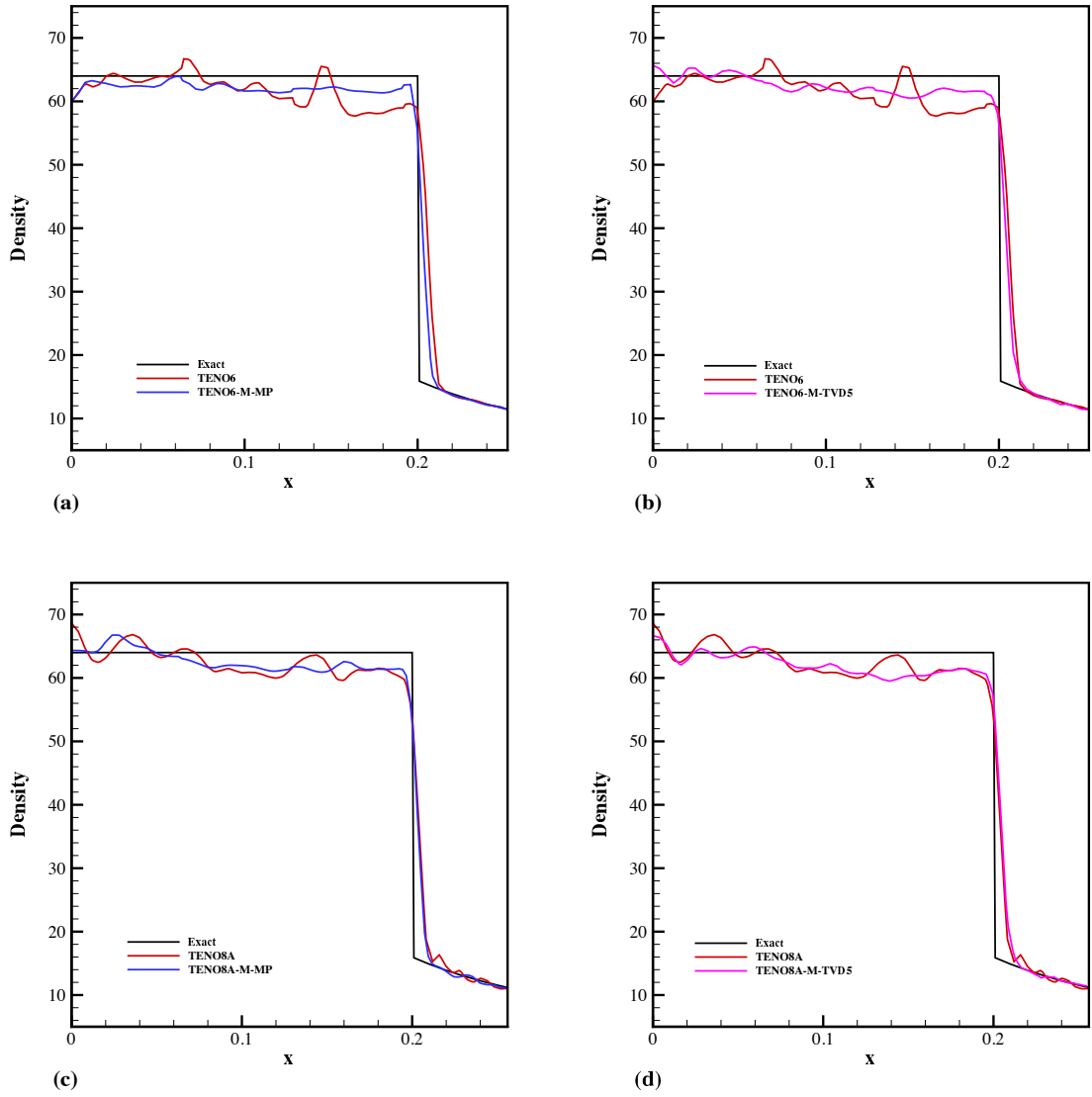}}
\caption{Noh problem: solutions from the TENO and TENO-M schemes. Top: density distribution (a, b) of six-point TENO6 and TENO6-M schemes; bottom: density distribution (c,d) of eight-point TENO8A and TENO8A-M schemes. Discretization is on $64\times64\times64$ uniformly distributed grid points and the final simulation time is $t = 0.6$.}
\label{Fig:Noh_1}
\end{figure}

As shown in Fig.~\ref{Fig:Noh_1}, the position of the shock wave is well captured and the density distributions from TENO and TENO-M show good agreement with the exact solution. Unlike the density drop caused by wall over-heating at the center observed with WENO schemes, e.g. see Fig.~10 in \cite{johnsen2010assessment}, no artificial wall heating appears with the TENO and TENO-M schemes. Moreover, compared with the classical TENO schemes, TENO-M-TVD5 and TENO-M-MP show an improved suppression of spurious oscillations.
}

{
\subsubsection{Le Blanc problem}

The Le Blanc problem \cite{loubere2005subcell} is a challenging case with very strong discontinues. The initial condition is
\begin{equation}
\label{eq:Leblanc}
(\rho ,u,p) = \left\{ {\begin{array}{*{20}{c}}
{(1,0,\frac{2}{3}\times10^{-1}),}&{\text{if } 0 \le x < 3} ,\\
{(10^{-3},0,\frac{2}{3}\times10^{-10}),}&{\text{if } 3 \le x \le 9} .
\end{array}} \right.
\end{equation}
The mesh resolution is $N = 900$ and the final simulation time is $t = 6$. The reference solutions are computed with the fifth-order WENO-JS scheme at a resolution $N = 2500$. The adiabatic coefficient is $\gamma = \frac{5}{3}$.

As shown in Fig.~\ref{Fig:Leblanc}, the TENO8A scheme leads to obvious numerical oscillations in density and velocity profiles. For the density profile, TENO8A-M-MP eliminates the overshoot at $x=6$ generated by the standard TENO8A scheme. In addition, TENO8A-M-MP preserves the monotonicity well in the velocity profile.
\begin{figure}[htbp]
\centering
\subfigure{\includegraphics[width=0.95\textwidth]{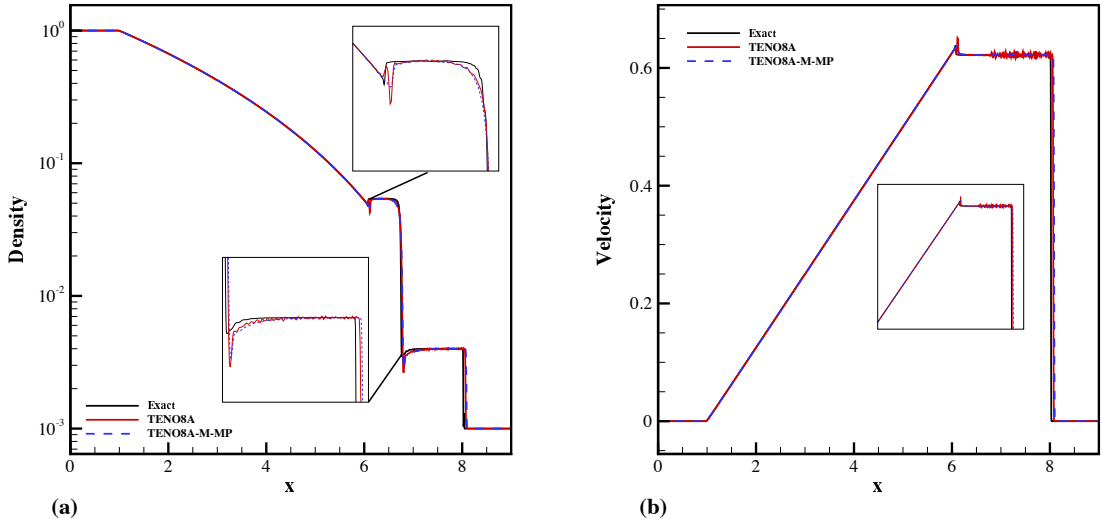}}
\caption{Le Blanc problem: solutions from the TENO and TENO-M schemes. Density distribution (a) and velocity distribution (b) of eight-point TENO8A and TENO8A-M schemes. Discretization is on $900$ uniformly distributed grid points and the final simulation time is $t = 6$.}
\label{Fig:Leblanc}
\end{figure}
}

\subsection{Computational efficiency}
In order to assess the computational efficiency of the proposed TENO-M scheme, we consider typical two-dimensional simulations, see Table.~\ref{Tab:efficiency}. And the cases involved are the Rayleigh-Taylor instability with a resolution of 64 $\times$ 256 and the double Mach reflection of a strong shock with a resolution of 800 $\times$ 200 respectively.  All simulations are conducted on the same desktop workstation. The results show that the cost of TENO8A-M scheme depends on the computational complexity of the limiters. For simple limiters, the computational efficiency of TENO8A-M can be even higher than that of TENO8A as the renormalization procedure of TENO is not needed.
The TENO8A-M-MP scheme is slightly more expensive than the TENO8A scheme by about $10\%$. For TVD-type limiters (TVD5 and VA), the computational efficiency of TENO8A-M is almost the same as the classical TENO8A scheme. Regarding to the good performance and flexibility of the new framework, we believe that the proposed method is promising for a wide range of applications.

\begin{table}[h]
  \centering
  \caption{The averaged computational time (in the parentheses) and the normalized values with respect to the computational time of TENO8A}
  \scriptsize
  \label{Tab:efficiency}
  \newcommand{\tabincell}[2]{\begin{tabular}{@{}#1@{}}#2\end{tabular}}
  \begin{tabular}{>{\centering\arraybackslash}m{2.8cm}
                  >{\centering\arraybackslash}m{2.0cm}
                  >{\centering\arraybackslash}m{1.5cm}
                  >{\centering\arraybackslash}m{2.5cm}
                  >{\centering\arraybackslash}m{2.5cm}
                  >{\centering\arraybackslash}m{1.8cm}}
  \toprule
  {Case} & {Grid number} &{TENO8A} &{TENO8A-MP}&{TENO8A-TVD5} &{TENO8A-VA}\\
  \hline

   {Rayleigh-Taylor instability} & $64\times256$ &  1(228.44s) &  1.088(248.58s) & 0.988(225.72s) & 0.997(227.89s)   \\

   {Double Mach reflection of a strong shock}  & $800\times200$ &  {1(1153.49s)} &  1.123(1295.31s) & 0.985(1137.11s) & 0.993(1145.91s)   \\

%   \hline
  \bottomrule

  \end{tabular}
  \end{table}

\section{Conclusions}
In this paper a flexible framework for constructing new shock-capturing TENO schemes is proposed. Six- and eight-point TENO-M schemes are developed, and their performance is demonstrated by conducting a set of critical benchmark cases. The conclusions are as follows.
\begin{itemize}
    \item The new framework establishes a unified concept of TENO schemes with classical nonlinear limiters for shock-capturing. Three stages are involved, (a) evaluating of candidate numerical fluxes and labelling each candidate stencil as smooth or nonsmooth by a ENO-like stencil selection procedure; (b) filtering the nonsmooth candidate stencils by an extra nonlinear limiter; (c) formulating the high-order reconstruction by combing the candidate stencils with optimal linear weights.

    \item  In smooth regions, all candidate stencils are identified as smooth by the TENO stencil selection procedure and consequently TENO-M recovers to TENO. In nonsmooth regions, stable shock-capturing capability is achieved since the nonsmooth candidates contributed to the final reconstructions are filtered to be oscillation-free by extra limiters.

    \item  The present framework can be applied to a wide range of limiters, such as TVD and MP. Such different nonlinear limiters are deployed straightforwardly to construct a new family of high-order TENO-M schemes. TENO-M schemes enable flexible control of dissipation in nonsmooth regions. A wide range of adjustable nonlinear dissipation can be considered in the current framework without affecting global dissipation and detriment to accuracy in low-wavenumber regions.

    \item  A set of critical benchmark cases is simulated. Numerical results demonstrate the capability of the new schemes in terms of recovering high-order accuracy in smooth regions, preserving low dissipation for resolution of fluctuations, and sharp capturing of discontinuities ({eliminating spurious oscillations particularly in extreme cases}). Moreover, different nonlinear limiters are tested. The results show that adapting nonlinear numerical dissipation in nonsmooth regions can be controlled by the choice of limiter function.

\end{itemize}
\section*{Acknowledgements}

The first author is supported by China Scholarship Council (NO. 201706290041). Lin Fu is funded by a CTR postdoctoral fellowship at Stanford University.
\appendix
\renewcommand{\appendixname}{Appendix}

\bibliographystyle{elsarticle-num-names}

\bibliography{explicit}

\end{document}